\documentclass[5p,times]{elsarticle}

\usepackage{lipsum}
\usepackage{amsmath}
\usepackage{amsfonts}
\usepackage{graphicx}
\usepackage{epstopdf}
\usepackage{algpseudocode}
\usepackage{algorithm}
\usepackage{algorithmicx}
\ifpdf%
  \DeclareGraphicsExtensions{.eps,.pdf,.png,.jpg}
\else
  \DeclareGraphicsExtensions{.eps}
\fi
\usepackage{amsopn}

\usepackage{booktabs}
\usepackage{tikz}
\usepackage{caption}
\usepackage{subcaption}
\usepackage{arydshln}
\usepackage{letltxmacro}
\usepackage{breqn}

\usepackage[justification=centering]{caption}

\newcommand{\eps}{\varepsilon}

\newcommand\ddfrac[2]{\frac{\displaystyle #1}{\displaystyle #2}}

\journal{Parallel Computing}

\newcommand{\ignore}[1]{}
\begin{document}
\begin{frontmatter}
\title{Low-Synch Gram-Schmidt with Delayed Reorthogonalization for Krylov Solvers}
%
\author{Daniel Bielich}
\author{Julien Langou}
\author{Stephen Thomas}
\author{Kasia \'{S}wirydowicz}
\author{Ichitaro Yamazaki}
\author{Erik G. Boman}

\begin{abstract}
The parallel strong-scaling of iterative methods is often 
determined by the number of global reductions at each iteration.
Low-synch Gram-Schmidt algorithms for nonsymmetric matrices are 
applied here to the GMRES and Krylov-Schur iterative solvers. 
The $QR$ decomposition is
``left-looking'' and processes one column at a time. 
Among the methods for generating an orthogonal basis
for the Krylov-Schur algorithm, the classical Gram Schmidt algorithm,
with reorthogonalization (CGS2) requires three global reductions per iteration.
A new variant of CGS2 that requires
only one reduction per iteration is applied to the Arnoldi-$QR$ iteration.
Delayed CGS2 (DCGS2) employs the minimum number of global reductions for a 
one-column at-a-time algorithm. The main idea behind the 
new algorithm is to group global reductions by rearranging the order of operations and thus DCGS2 must be carefully integrated into an Arnoldi--$QR$ 
expansion. Numerical stability experiments assess robustness for 
Krylov-Schur eigenvalue computations.
Performance experiments on the ORNL Summit supercomputer then 
establish the superiority of DCGS2.
\end{abstract}

\end{frontmatter}
\section{Introduction}\label{db:sec:intro}

\begin{center}
The Rolling Stones \cite{rollingStones} - You Can't Always Get What You Want ...\\
~\\
No, you can't always get what you want\\
You can't always get what you want\\
You can't always get what you want\\
But if you try sometimes, well you just might find\\
You get what you need\\
\end{center}

Let $A$ be an $m \times m$ real-valued matrix.
In this manuscript, $A$ is employed in two parallel
computations: (1) finding a solution of the linear system $Ax=b$ with a Krylov
subspace method such as GMRES \cite{Saad86} and (2) finding the eigenvalues of
$A$ using Krylov-Schur \cite{Stewart86}.  In both instances, an orthogonal
basis for the Krylov subspace ${\cal K}_n$ is required. The size of the basis
is $n \ll m$. 

Gram-Schmidt produces a $QR$ decomposition of a matrix $A$ and for
Arnoldi--$QR$, the factorized matrix is $B = [\:r_0,\:AQ_n\:]$. The Arnoldi
algorithm applies the Gram-Schmidt process to produce a linearly independent
basis consisting of the columns of $Q_n$ in the Arnoldi expansion
$Q_n^TAQ_n=H_n$.  Krylov methods for linear system and eigenvalue solvers both
depend upon the orthogonality of the basis for the Krylov subspace ${\cal K}_n(B)$
in finite precision arithmetic.  

The loss of orthogonality of the computed basis -- as measured by $\|I -
Q_n^TQ_n\|_F$ -- may deviate substantially from machine precision ${\cal
O}(\eps)$, (see Giraud et al.~\cite{gilr:cma:05}).  
When linear independence is completely lost, the Krylov iterations,
may fail to converge.  For the solution of linear systems of equations $Ax =
b$, Paige et al.~\cite{2006--simax--paige-rozloznik-strakos} show that when the
loss of orthogonality is ${\cal O}(\eps)\kappa(B)$, then MGS-GMRES is backward
stable for the solution of linear systems. Here, $\kappa(B)$ is the condition
number $\kappa(B) = \sigma_{\max}(B)/\sigma_{\min}(B)$, where $\sigma_i(B)$ are
the singular values of the matrix $B$. For eigenvalue computations employing
Krylov-Schur, accurate and converged eigenvalue approximations are obtained
when the loss of orthogonality remains close to ${\cal O}(\eps)$.  In this
paper, a stable Arnoldi algorithm is presented that performs well on an
Exascale class supercomputer.

Krylov linear system and eigenvalue solvers are often required for extreme
scale physics simulations and implemented on parallel (distributed memory)
machines.  Their strong-scaling is limited by the number and frequency of
global reductions, in the form of {\tt MPI\_AllReduce}. These communication
patterns are expensive \cite{BienzGroppOlson2019}.   Our new algorithms are
designed such that they require only one reduction to normalize each vector and
apply projections.  The focus here is on methods that process the Krylov
vectors one column at a time as opposed to blocks (e.g. only one column becomes
available at a time, and thus the vectors are processed in a ``left-looking''
fashion).  As indicated previously,  the Krylov basis size $n$, is assumed to
be much smaller than the dimension $m$ of the matrix $A$, one can think of $m$
as infinite and computations such as inner-products are avoided as often as
possible. $Q_n$ is then referred to as a tall and skinny matrix. These are
typical of Krylov iterative methods for solving a linear system of equations or
eigenvalue problems, which rely on the Arnoldi expansion.

Classical Gram-Schmidt (CGS) is preferred for a single program multiple data
(SPMD) model of computation because it requires only two global reductions for
each column vector (a projection step, followed by a vector normalization).  In
practice, however, CGS leads to numerical instability for the solution of
$Ax=b$ and also eigenvalues, because the loss of orthogonality is
$\mathcal{O}(\eps)\kappa^2(A)$.  This bound was conjectured for a long time and
finally proven in two papers \cite{glre:nm:05,smbl:nm:06}.  The GMRES iteration
will stall and fail to converge if linear independence of the Krylov vectors
is completely lost, for example, when $\|S\|_2 = 1$ as described by Paige
\cite{paige2018-Book-}, where the matrix $S$ was introduced in Paige et
al.~\cite{2006--simax--paige-rozloznik-strakos}.  In order to obtain backward
stable eigenvalues from Krylov-Schur,  Stewart \cite{Stewart86} demonstrates
that $\mathcal{O}(\eps)$ loss of orthogonality suffices.

To reduce the loss of the orthogonality to machine precision ${\cal O}(\eps)$,
the CGS algorithm can be applied twice (CGS2) to reorthogonalize the basis
vectors.  This is the ``twice is enough'' result from Kahan and Parlett
\cite{doi:10.1137/1.9781611971163}, which has also been proven by Giraud et
al.~\cite{glre:nm:05}.  Given the assumption that $c\eps\kappa(A) < 1$ for a
given $m\times n$ input matrix $A$ and constant $c=\mathcal{O}(m^2n^3)$, then
CGS2 can construct orthogonal columns to machine precision, Theorem 2 in
\cite{glre:nm:05}. The number of floating point operations (flops) for CGS2 is
therefore $4mn^2$ (twice the cost of CGS), and requires three global
reductions.

A one-reduce variant of CGS2 was derived in \cite{2020-swirydowicz-nlawa} and
is applied to the Arnoldi-$QR$ iteration in this paper. It achieves the same
$\mathcal{O}(\eps)$ loss of orthogonality as the original CGS2 but requires
only one reduction per column vector.  To minimize the number of global
reductions and avoid cancellation errors, the normalization step is lagged and
the Pythagorean theorem is employed.  The reorthogonalization is also delayed
to the next iteration and thus is performed ``on-the-fly" as opposed to a
second pass.  The resulting algorithm combines these two steps into one global
reduction and is referred to as the delayed DCGS2 (this is explained in detail
in Section 3.).

Extensive numerical results are presented for the Krylov-Schur eigenvalues to
demonstrate the numerical stability and accuracy of DCGS2-Arnoldi.
Strong-scaling results are presented for the ORNL Summit supercomputer to
demonstrate that the DCGS2 algorithm improves the CGS2 compute times by a
factor of up to $2\times$ on many-core architectures such as GPUs, while 
maintaining the same loss of orthogonality as the original CGS2-Arnoldi algorithm.

\section{Low-Synch Gram-Schmidt Algorithms}
The development of low-synch MGS and CGS2 was largely driven by applications
that need stable, yet scalable solvers. Even though MGS-GMRES is backward
stable for the solution of linear systems, CGS2-GMRES was found to be more
scalable for massively parallel computation on the Cray T3D in a 1998 study by
Frayss\'e et al.~\cite{Frayse98} and included in the Trilinos framework by
Bavier et al.~\cite{Belos}.  The more recent development of a one-reduce
MGS-GMRES by Swirydowicz et al.~\cite{2020-swirydowicz-nlawa} implies that a
re-evaluation of these results is certainly warranted in the context of a
one-reduce DCGS2. 

As already stated, the CGS2 algorithm requires three global reductions per
iteration: one for the first projection, another for the second pass and a
third for the normalization. The one-reduce DCGS2 delays reorthogonalization.
This is achieved by lagging the normalization as originally proposed by  Kim
and Chronopoulos \cite{Kim92}) and then applying Stephen's trick. The Pythagorean trick introduced by Smoktunowicz et al.~\cite{smbl:nm:06} avoids cancellation errors and Carson et al.~\cite{2020_arXive_Carson} generalize this to block Gram-Schmidt algorithms.

The delayed normalization for the Arnoldi iteration was employed by Hernandez
et al.~\cite{2007-hernandez-parco} (Arnoldi with Delayed Reorthogonalization~-~ADR) 
without a correction after normalization and we refer to their Gram-Schmidt
algorithm as DCGS2-HRT. DCGS2-Arnoldi as derived by these authors
is not forward stable for eigenvalue
computations because the loss of orthogonality is at least ${\cal
O}(\eps)\kappa(B)$, see Section \ref{num}.  Delaying the normalization also
requires a scaling to update the Krylov vectors in the Arnoldi expansion.
Recent work
\cite{2020-swirydowicz-nlawa,2020-yamazaki-proceedings-of-siam-pp20} describes a
one-reduce inverse compact $WY$ MGS algorithm with a triangular solve in the
projection step.  This ICWY--MGS requires only one reduction per iteration,
the same as DCGS2.  The original modified Gram-Schmidt (MGS) requires $2mn^2$
flops (the same as CGS) and applies the elementary rank-1 projections $I -
q_jq_j^T$ sequentially, requiring a separate global reduction for each
inner-product. The complexity of ICWY-MGS is also $3mn^2$, including an additional 
$mn^2$ flops for constructing the triangular factor $L$. The loss of
orthogonality for MGS is  $\mathcal{O}(\eps)\kappa(A)$ (see Bj\"orck
\cite{Bjork_1967}).  Whereas, it is $\mathcal{O}(\eps)$ for CGS2. 

An MGS2, analogous to CGS2, exists.  In practice, MGS2 exhibits an
$\mathcal{O}(\eps)$ loss of orthogonality.  The number of flops for MGS2 is
$4mn^2$.  (Double the cost of CGS or MGS and the same as CGS2).  Low synch MGS2
requires $4 mn^2$ flops ($2mn^2$ for the first pass and $2mn^2$ for the second
pass).  These costs will be important considerations in strong scaling studies
of these new algorithms on the ORNL Summit supercomputer.

\section{DCGS2 Algorithm for the $QR$ Decomposition}
In this section, the classical Gram-Schmidt algorithm to compute the $QR$
decomposition of an $m\times n$ matrix $A$ is reviewed.
Algorithm \ref{db:alg:cgs2} displays the steps for the $j$--th iteration of
CGS2. 
The column vector ${a}_j$ is twice orthogonally projected onto the orthogonal complement of $Q_{1:j-1}$, then normalized
\begin{eqnarray} 
\nonumber {w}_{j} & = & \left(\:I - Q_{1:j-1}\:Q_{1:j-1}^T \:\right)\:{a}_j\\
\nonumber
& = & a_j - Q_{1:j-1}\:S_{1:j-1, j} \quad\mbox{ where }\quad 
S_{1:j-1, j} = Q_{1:j-1}^T  {a}_j,
\end{eqnarray}
followed by a second application of the projector in the form 
\begin{eqnarray}
\nonumber
  {u}_{j} & = & w_j - Q_{1:j-1}\:C_{1:j-1, j} \quad\mbox{ where }
  \quad C_{1:j-1, j} = Q_{1:j-1}^T  {w}_j.
\end{eqnarray}
Finally, the vector $u_j$ is normalized to produce the vector $q_j$,
\[
  {q}_j = {u_j}/\alpha_j \quad\mbox{ where }\quad \alpha_j = \|\:u_j\:\|_2.
\]
The $QR$ decomposition of $A_{1:j} = Q_{1:j} R_{1:j, 1:j}$, 
is produced at the end of the $j$--th iteration, where
\[
  R_{1:j-1, j} = S_{1:j-1, j} + C_{1:j-1, j}
  \quad\mbox{ and }\quad
  R_{j,j} = \alpha_j.
\]
The $j$--th iteration of DCGS2  is displayed in Algorithm~\ref{db:alg:cgs2}.
The three global
reductions appear in steps 1, 3, and 5. 

\begin{algorithm}[htb]
\begin{algorithmic}[1]
\Statex{\hspace{\algorithmicindent} \emph{// first projection}}
\State{\hspace{\algorithmicindent} $S_{1:j-1,j} =  {Q}_{1:j-1}^T {a}_{j}$       \hspace{8mm} \emph{// global reduction}}
\State{\hspace{\algorithmicindent} ${w}_j =  {a}_j -  {Q}_{1:j-1} {S}_{1:j-1,j}$}
\Statex
\Statex{\hspace{\algorithmicindent} \emph{// second projection}}
\State{\hspace{\algorithmicindent} $C_{1:j-1,j} =  {Q}_{1:j-1}^T {w}_{j}$       \hspace{7mm} \emph{// global reduction}}
\State{\hspace{\algorithmicindent} ${u}_{j} =  {w}_{j} -  {Q}_{1:j-1} {C}_{1:j-1,j}$}
\Statex
\Statex{\hspace{\algorithmicindent} \emph{// normalization}}
\State{\hspace{\algorithmicindent} $\alpha_j = \|\:u_j\:\|_2$   \hspace{16.4mm} \emph{// global reduction}}
\State{\hspace{\algorithmicindent} ${q}_{j} = {u}_{j} / \alpha_j$}
\Statex
\Statex{\hspace{\algorithmicindent} \emph{// representation $ R_j$}}
\State{\hspace{\algorithmicindent} $R_{1:j-1,j} =  {S}_{1:j-1,j} +  {C}_{1:j-1,j}$}
\State{\hspace{\algorithmicindent} $R_{j,j} = \alpha_j$ }
\end{algorithmic}
\caption{\label{db:alg:cgs2} Classical Gram-Schmidt with reorthogonalization (CGS2) }
\end{algorithm}

The $j$--th iteration of DCGS2  is displayed in Algorithm~\ref{db:alg:qr}.  In
order to perform one global reduction,  the second projection and normalization
are lagged or delayed to the next iteration.  The purpose of the Pythagorean
trick is to mitigate cancellation errors due to finite precision arithmetic.
Namely, the norm of the updated vector $u_j$ is computed as follows
\begin{eqnarray}
\nonumber
\alpha_j^2 
& = & \left(\:{w}_{j} -  {Q}_{1:j-1}\:{C}_{1:j-1,j}\:\right)^T \left(\:{w}_{j} -  {Q}_{1:j-1}\:{C}_{1:j-1,j}\:\right)\\
\nonumber
& = & {w}_{j}^T {w}_{j}  -  2 {C}_{1:j-1,j}^T\:{C}_{1:j-1,j} \\ \nonumber 
& + & {C}_{1:j-1,j}^T\: \left(\:{Q}_{1:j-1}^T\:{Q}_{1:j-1}\:\right)\: {C}_{1:j-1,j}\\
\nonumber
& = & {w}_{j}^T {w}_{j}  -  {C}_{1:j-1,j}^T\:{C}_{1:j-1,j}
\end{eqnarray}
where $C_{1:j-1,j} =  {Q}_{1:j-1}^T\:{w}_{j}$ and the orthogonality of
$Q_{1:j-1}$ is assumed in finite precision arithmetic to be ${\cal O}(\eps)$.
This corresponds to Step 3 in Algorithm \ref{db:alg:qr}.

Because the normalization is delayed, the Pythagorean trick allows us to
compute the norm by utilizing ${w}_{j-1}$ instead of ${u}_{j-1}$.
For the column vector ${a}_j$, the scalar $S_{j-1,j}$ is
computed with $S_{j-1, j} =  {w}_{j-1}^T\:{a}_j$ instead of
${q}_{j-1}^T\:{a}_j$.  This is Stephen's trick, it captures the computation
${q}_{j-1}^T\:{a}_j$ using ${w}_{j-1}$ instead of $q_{j-1}$ and results
in a corrected projection step within the Gram Schmidt process, 
\begin{eqnarray}
q_{j-1}^T\:a_j
\nonumber  & =  & \frac{1}{\alpha_{j-1}}\left(\: {w}_{j-1} - Q_{1:j-2} {C}_{1:j-2,j-1} \:\right)^T  {a}_j  \\
\nonumber  & = & \frac{1}{\alpha_{j-1}}
\left(\: {w}_{j-1}^T  {a}_j -  {C}_{1:j-2,j-1}^T\:Q_{1:j-2}^T\:{a}_j \:\right) \\
\nonumber  & = & \frac{1}{\alpha_{j-1}}
\left(\: S_{j-1,j} - {C}_{1:j-2,j-1}^T\:{S}_{1:j-2,j}  \:\right) 
\end{eqnarray}
The correction of the vector norm corresponds to Steps 3 and 5
in Algorithm \ref{db:alg:qr} given below.
\begin{algorithm}[htb]
\begin{algorithmic}[1]
\State{\hspace{\algorithmicindent}  
$[\: {Q}_{1:j-2},\:  {w}_{j-1}\:]^T\:[\: {w}_{j-1},\:  {a}_{j} \:]$
 \hspace{7mm} \emph{// global reduction}}
\Statex
\Statex{\hspace{\algorithmicindent} ${S}_{1:j-2,j} =  {Q}_{1:j-2}^T\:{a}_{j}$ \hspace{5mm} \quad and \quad $S_{j-1,j} =  {w}_{j-1}^T {a}_{j}$ }
\Statex{\hspace{\algorithmicindent} ${C}_{1:j-2,j-1} =  {Q}_{1:j-2}^T\:{w}_{j-1}$ \quad and \quad  $\beta_{j-1} =  {w}_{j-1}^T {w}_{j-1}$ }
\Statex
\Statex{\hspace{\algorithmicindent} \emph{// delayed reorthogonalization}}
\State{\hspace{\algorithmicindent} ${u}_{j-1} =  {w}_{j-1} -  {Q}_{1:j-2}\:{C}_{1:j-2,j-1}$ }
\Statex
\Statex{\hspace{\algorithmicindent} \emph{// delayed normalization}}
\State{\hspace{\algorithmicindent}  $\alpha_{j-1} = \left\{\: \beta_{j-1} - {C}^T_{1:j-2,j-1}\:{C}_{1:j-2,j-1}\:\right\}^{1/2}$ }
\State{\hspace{\algorithmicindent} ${q}_{j-1} = \ddfrac{1}{\alpha_{j-1}}\:{u}_{j-1} $ }
\Statex
\Statex{\hspace{\algorithmicindent} \emph{// projection}}
\State{\hspace{\algorithmicindent} $S_{j-1,j} = \ddfrac{1}{\alpha_{j-1}}
\left(\: S_{j-1,j} -  {C}^T_{1:j-2,j-1}\:S_{1:j-2,j}\:\right) $ }
\State{\hspace{\algorithmicindent} ${w}_j =  {a}_j -  {Q}_{1:j-1}\:{S}_{1:j-1,j}$ 
\Statex
\Statex{\hspace{\algorithmicindent} \emph{// representation $ R_{j-1}$}}
\State{\hspace{\algorithmicindent} $R_{1:j-2,j-1} =  {S}_{1:j-2,j-1} +  {C}_{1:j-2,j-1}$ }
\State{\hspace{\algorithmicindent} $R_{j-1,j-1} = \alpha_{j-1}$ }}
\end{algorithmic}
\caption{\label{db:alg:qr} Delayed Classical Gram-Schmidt with reorthogonalization (DCGS2)}
\end{algorithm}
For the $n$--th iteration, CGS2 is applied
and incurs two additional global reductions.
\section{DCGS2 Algorithm for the Arnoldi Expansion}
\label{db:sec:arnoldi}
Algorithm \ref{db:alg:cgs2_arnoldi} displays CGS2 for the Arnoldi expansion.
The only difference from the $QR$ decomposition in Algorithm \ref{db:alg:cgs2}
is that the next basis vector ${v}_j$ is generated by applying a matrix-vector
product to the previously normalized column vector ${q}_{j-1}$.  At the end of
iteration $j$$-1$, in exact arithmetic, the matrices would satisfy the Arnoldi
expansion,
\begin{equation}\label{eq:arnoldi}
A\: Q_{1:j-2} = Q_{1:j-1}\: H_{1:j-1,1:j-2}.
\end{equation}
\begin{algorithm}[htb]
\begin{algorithmic}[1]
\Statex{\hspace{\algorithmicindent} \emph{// generation of next vector}}
\State{\hspace{\algorithmicindent} ${v}_{j} = A\:{q}_{j-1}$}
\Statex
\Statex{\hspace{\algorithmicindent} \emph{// first projection}}
\State{\hspace{\algorithmicindent} ${S}_{1:j-1,j} =  {Q}_{1:j-1}^T\:{v}_{j}$ \hspace{8.2mm} \emph{// global reduction}}
\State{\hspace{\algorithmicindent} ${w}_j =  {v}_j -  {Q}_{1:j-1}\:{S}_{1:j-1,j}$}
\Statex
\Statex{\hspace{\algorithmicindent} \emph{// second projection}}
\State{\hspace{\algorithmicindent} ${C}_{1:j-1,j} =  {Q}_{1:j-1}^T\:{w}_{j}$ \hspace{7mm} \emph{// global reduction}}
\State{\hspace{\algorithmicindent} ${u}_{j} =  {w}_{j} -  {Q}_{1:j-1}\:{C}_{1:j-1,j}$}
\Statex
\Statex{\hspace{\algorithmicindent} \emph{// normalization}}
\State{\hspace{\algorithmicindent} $\alpha_j = \|\:{u}_{j}\:\|_2$\ \hspace{16.8mm} \emph{// global reduction}}
\State{\hspace{\algorithmicindent} ${q}_{j} = \ddfrac{1}{\:\alpha_j}\: {u}_{j}$}
\Statex
\Statex{\hspace{\algorithmicindent} \emph{// representation $H_j$}}
\State{\hspace{\algorithmicindent} ${H}_{1:j-1,j} =  {S}_{1:j-1,j} +  {C}_{1:j-1,j}$}
\State{\hspace{\algorithmicindent} $H_{j,j} = \alpha_j$ }
\end{algorithmic}
\caption{\label{db:alg:cgs2_arnoldi} Arnoldi-$QR$ (CGS2)}
\end{algorithm}
A one-reduction DCGS2-Arnoldi will now be derived.  The representation error
and loss of orthogonality are maintained at the same level as the
CGS2-Arnoldi.  

With lagged vector updates, the next basis vector is generated by applying a
matrix-vector product to the current vector.  Namely, the next vector ${v}_j$
is computed as $A\:{w}_{j-1}$ by using the vector $ {w}_{j-1}$ instead of
${q}_{j-1}$, where $ {q}_{j-1}$ is the previously constructed orthogonal
column.  Thus, an effective strategy is required to compute $w_{j}$ from
$Aw_{j-1}$ and also to generate the Hessenberg matrix $H_j$ in the Arnoldi
expansion.

After a delay of one iteration, the vector $q_{j-1}$, is computed using
$w_{j-1}$ as follows
\begin{eqnarray}
\label{08}  &&  {q}_{j-1} = \frac{1}{\alpha_{j-1}}\left(\:{w}_{j-1} -  {Q}_{1:j-2}  {C}_{1:j-2,j-1}\:\right)
\end{eqnarray}
Equation \eqref{08} may also be interpreted as a $QR$ factorization of
the matrix ${W}_{1:j-1}$, with columns $[\:w_1,\: \ldots, \:w_{j-1} \:]$
\begin{equation}
\label{09-matrix}  {Q}_{1:j-1} =  {W}_{1:j-1} C_{1:j-1,1:j-1}^{-1},
\end{equation}
where $C$ is an upper triangular matrix.

Multiplying \eqref{08} by ${A}$ from the left, it follows that
\begin{eqnarray}
{v}_{j}   & = &  {A}\:{q}_{j-1}\\
\nonumber    
& = & \frac{1}{\alpha_{j-1}}\left(\:{A}\:{w}_{j-1} -  {A}\:{Q}_{1:j-2}\:{C}_{1:j-2,j-1}\:\right) \\
\label{11} 
& = & \frac{1}{\alpha_{j-1}}\left(\:{A}\:{w}_{j-1} -  {Q}_{1:j-1}\:H_{1:j-1,1:j-2}\:{C}_{1:j-2,j-1}\:\right) 
\nonumber    
\end{eqnarray}
Next the vector ${w}_j$ is computed, which is the vector produced after
projection of ${v}_j$ onto the basis vectors in ${Q}_{1:j-1}$,
\begin{eqnarray}
\label{06a} 
{w}_j & = &  {A}\:{q}_{j-1}  -  {Q}_{1:j-1}\:{Q}_{1:j-1}^T\:{A}\:{q}_{j-1}\\
\nonumber    
& = & \frac{1}{\alpha_{j-1}} \left(\: {A} {w}_{j-1} -  {Q}_{1:j-1} H_{1:j-1,1:j-2}  {C}_{1:j-2,j-1} \:\right) \\
\nonumber    
&   & -  {Q}_{1:j-1}  {Q}_{1:j-1}^T \\
\nonumber 
& \times & \:\frac{1}{\alpha_{j-1}}
\left(\:{A} {w}_{j-1} -  {Q}_{1:j-1}\:H_{1:j-1,1:j-2}\:{C}_{1:j-2,j-1}\: \right)\\
\nonumber    
w_j & = &  
\frac{1}{\alpha_{j-1}}\left(\:{A}\:{w}_{j-1}  -  {Q}_{1:j-1}\:{Q}_{1:j-1}^T\:{A}\:{w}_{j-1}\:\right) \\
& + & \frac{1}{\alpha_{j-1}}\: Q_{1:j-1}\:
\left(\:I -  {Q}_{1:j-1}^T\: {Q}_{1:j-1}\:\right)\:H_{1:j-1,1:j-2}\:{C}_{1:j-2,j-1} 
\nonumber
\end{eqnarray}
The last term is dropped from \eqref{06a},
%
%
for two reasons,
\begin{itemize}
\item DCGS2 is constructed such that the loss of orthogonality 
$\|I -  {Q}_{1:j-1}^T\:{Q}_{1:j-1}\|_F$ is ${\cal O}(\eps)$, and
\item ${C}_{1:j-2,j-1}/\alpha_{j-1}$ 
is expected to be ${\cal O}(\eps)\kappa(A)$. Hence, when 
$\kappa(A) \le {\cal O}(1/\eps)$, the norm of the term is ${\cal O}(1)$.
\end{itemize}
Therefore, at this point (\ref{06a}) becomes an approximation and
\begin{eqnarray}
\nonumber   
{w}_j & = & \frac{1}{\alpha_{j-1}}
\left(\:{A}\:{w}_{j-1}  - {Q}_{1:j-1}\:{Q}_{1:j-1}^T\:{A}\:{w}_{j-1}\:\right) \\
\nonumber
& = & \frac{1}{\alpha_{j-1}}
\left(\:{A}\:{w}_{j-1}  - {Q}_{1:j-2}\:{Q}_{1:j-2}^T\:{A}\:{w}_{j-1}\:\right) \\
\nonumber         
& - & \frac{1}{\alpha_{j-1}}\:{q}_{j-1}\:{q}_{j-1}^T\:{A}\:{w}_{j-1}
\label{31}        
\nonumber
\end{eqnarray}
Noting that ${S}_{1:j-2,j} =  {Q}_{1:j-2}^T\:{A}\:{w}_{j-1}$, it follows that.
\begin{eqnarray}
\nonumber 
{w}_j & = &
\frac{1}{\alpha_{j-1}}
\left(\:{A}\:{w}_{j-1} - {Q}_{1:j-2}\:{S}_{1:j-2,j} -  {q}_{j-1}\:{q}_{j-1}^T\:{A}\:{w}_{j-1}\:\right) 
\nonumber
\end{eqnarray}
Finally, from \eqref{08}, it is possible to compute 
%
\begin{eqnarray}
\nonumber  
{q}_{j-1}^T\:{A}\:{w}_{j-1} 
& = & \frac{1}{\alpha_{j-1}}
\left(\:{w}_{j-1} -  {Q}_{1:j-2}\:{C}_{1:j-2,j-1}\:\right) ^T  {A} {w}_{j-1}  \\
\nonumber
& = & \frac{1}{\alpha_{j-1}}
\left(\:{w}_{j-1}^T {A}\:{w}_{j-1} - 
{C}_{1:j-2,j-1}^T\:{Q}_{1:j-2}^T\:{A}\:{w}_{j-1}\: \right)   \\
\label{30}
& = & \frac{1}{\alpha_{j-1}}
\left(\: S_{j-1,j} - {C}_{1:j-2,j-1}^T\:{S}_{1:j-2,j} \:\right)
\nonumber
\end{eqnarray}
where $S_{j-1,j} =  {w}_{j-1}^T {A} {w}_{j-1}$.  
This step is Stephen's trick in the context of Arnoldi.

After substitution of this expression, it follows that
\begin{eqnarray}
\label{32}   
{w}_j & = & \frac{1}{\alpha_{j-1}}
\left(\:{A}\:{w}_{j-1}  - {Q}_{1:j-2}\:{S}_{1:j-2,j}\:\right) \\
\nonumber
& - & \frac{1}{\alpha_{j-1}^2}
\:{q}_{j-1}\:\left(\: S_{j-1,j} -  {C}_{1:j-2,j-1}^T\:{S}_{1:j-2,j}\:\right) \\
\nonumber
& = & \frac{1}{\alpha_{j-1}}\:{A}\:{w}_{j-1} - {Q}_{1:j-1}\:{T}_{1:j-1,j}
\nonumber
\end{eqnarray}
where
\[
{T}_{1:j-2,j} = \frac{1}{\alpha_{j-1}} \:{S}_{1:j-2,j} 
\]
and
\[
T_{j-1,j} = \frac{1}{\alpha_{j-1}^2}\:\left(\: S_{j-1,j} -  {C}_{1:j-2,j-1}^T\:{S}_{1:j-2,j}\:\right).
\]

The $(j-1)$--th column of the Hessenburg matrix $H$ is computed as follows and
satisfies the Arnoldi relation \eqref{eq:arnoldi}. First, reorder \eqref{32} 
into a factorization form
\begin{eqnarray}
\label{41}
{A}\:{w}_{j-1} & = &  {Q}_{1:j-2}\:{S}_{1:j-2,j}  \\
\nonumber 
& + & 
\frac{1}{\alpha_{j-1}}\:{q}_{j-1} \:
\left(\:S_{j-1,j} -   {C}_{1:j-2,j-1}^T\:{S}_{1:j-2,j}\:\right) + \alpha_{j-1}\:{w}_j 
\end{eqnarray}
%
%
From (\ref{08}), it also follows that
\begin{eqnarray}
\label{09} {w}_{j} 
& = &  {Q}_{1:j-1}\:{C}_{1:j-1,j} +   \alpha_j\:{q}_{j}
\end{eqnarray}
which represents the orthogonalization of the vector ${w}_{j}$.
By replacing $ {w}_j$ in \eqref{41} with the expression in \eqref{09}, obtain
\begin{eqnarray}
\nonumber   
{A} {w}_{j-1} 
& = & 
{Q}_{1:j-2}\:{S}_{1:j-2,j} \\
\nonumber  
& + & 
\frac{1}{\alpha_{j-1}} \:{q}_{j-1} \:
\left(\: S_{j-1,j} -  {C}_{1:j-2,j-1}^T  {S}_{1:j-2,j} \: \right) \\
\nonumber  
& + & \alpha_{j-1}\:{Q}_{1:j-1}\:{C}_{1:j-1,j} +  \alpha_j\:\alpha_{j-1}\:{q}_{j}\: \\
\label{Aw} 
{A} {w}_{j-1}
& = & {Q}_{1:j-2} \left(\:  {S}_{1:j-2,j} +  \:\alpha_{j-1}\:{C}_{1:j-2,j} \:\right) \\
\nonumber  
& + &  \frac{1}{\alpha_{j-1}} \: {q}_{j-1} \:
\left(\:  S_{j-1,j} -  {C}_{1:j-2,j-1}^T  {S}_{1:j-2,j}\: \right) \\
\nonumber  & + &  \alpha_{j-1} \: C_{j-1,j} \:{q}_{j-1}    + 
\alpha_j\:\alpha_{j-1}\:{q}_{j} .
\nonumber
\end{eqnarray}
This is the representation of ${A}\,{w}_{j-1}$ in the Krylov subspace spanned
by the orthogonal basis vectors ${Q}_{1:j}$ using the matrices $S$ and $C$.
However, the representation of ${A}\:{q}_{j-1}$ in $ {Q}_{1:j}$ 
with the matrix $H$ is still required.  Namely, write \eqref{11} as
\begin{eqnarray}
\nonumber
{A}  {q}_{j-1} & = & \ddfrac{1}{\alpha_{j-1}} 
\left(\:{A} {w}_{j-1}  - {Q}_{1:j-2}\: H_{1:j-2,1:j-2} \: {C}_{1:j-2,j-1} \:\right) \\
\nonumber 
& - & \frac{1}{\alpha_{j-1}}\:{q}_{j-1}\: H_{j-1,j-2}\: C_{j-2,j-1}
\label{Aq}
\end{eqnarray}
$H_{j-1,j-1}$ is now computed using ${C}_{j-2}$ and $H_{1:j-1,j-2}$.

Replacing ${A} {w}_{j-1}$ with \eqref{Aw}, it follows that
%
%
\begin{eqnarray}
\label{52}   
{A} {q}_{j-1} 
& = & 
{Q}_{1:j-2}\:\left(\: \frac{1}{\alpha_{j-1}}\: S_{1:j-2,j} + C_{1:j-2,j} \:\right) \\
\nonumber  
& + &  \frac{1}{\alpha_{j-1}^2} \: {q}_{j-1} \:
\left(\: S_{j-1,j} - C_{1:j-2,j-1}^T\:S_{1:j-2,j}\: \right) \\
\nonumber 
& + & \alpha_{j}\:{q}_{j} - \frac{1}{\alpha_{j-1}}\:{Q}_{1:j-2}\: H_{1:j-2,1:j-2}\:C_{1:j-2,j-1}  \\
\nonumber  
& + & C_{j-1,j}\:{q}_{j-1}\: - \frac{1}{\alpha_{j-1}}\: 
\:H_{j-1,j-2}\:C_{j-2,j-1} \: {q}_{j-1}
\end{eqnarray}
To summarize
\begin{enumerate}
\item ${A}\:{q}_{j-1} =  {Q}_{1:j-2}\:{A}_{1:j-2,j} +  
{q}_{j-1}\: s_{j-1,j} +  \alpha_j\: {q}_{j} $ 
is a standard $QR$ decomposition obtained by a Gram-Schmidt process.
\item ${C}_{1:j-2,j}$ and $C_{j-1,j}$ are the standard 
reorthogonalization terms in the representation equation,
\item ${C}_{1:j-2,j-1}^T\:{S}_{1:j-2,j}$  is Stephen's trick.
\item $H_{1:j-2,1:j-2}\:{C}_{1:j-2,j-1}$ and $H_{j-1,j-2}\:C_{j-2,j-1}$ are the
representation error correction terms.
\item $\ddfrac{1}{\alpha_{j-1}}$ is due to using unnormalized quantities
and these must be corrected by scaling.
\end{enumerate}
Items 1, 2 and 3 are present in both the $QR$ decomposition and Arnoldi expansion. 
Items 4 and 5 are specific to Arnoldi. 

According to (\ref{52}), in order to obtain the $(j-1)$--th column of the
Arnoldi relation \eqref{eq:arnoldi}, the column $H_{1:j,j-1}$ is computed as follows
\begin{eqnarray}
\nonumber   
{H}_{1:j-2,j-1} & = & 
\ddfrac{1}{\alpha_{j-1}}\:{S}_{1:j-2,j} +  {C}_{1:j-2,j} \\
\nonumber 
& - & \ddfrac{1}{\alpha_{j-1}}\: H_{1:j-2,1:j-2}\:{C}_{1:j-2,j-1}\\
\nonumber
& = &  {T}_{1:j-2,j} +  {C}_{1:j-2,j} - 
\ddfrac{1}{\alpha_{j-1}}\: H_{1:j-2,1:j-2}\: {C}_{1:j-2,j-1} \\
%
%
\nonumber  
H_{j-1,j-1} & = & 
\ddfrac{1}{\alpha_{j-1}^2} \:
\left(\: S_{j-1,j} -  {C}_{1:j-2,j-1}^T\:{S}_{1:j-2,j}\: \right) + C_{j-1,j} \\
\nonumber 
& - & \ddfrac{1}{\alpha_{j-1}}\: H_{j-1,j-2}\:C_{j-2,j-1}\\
\nonumber  
H_{j,j-1}  & = & \alpha_j
\end{eqnarray}
Finally, the DCGS2--Arnoldi is presented in Algorithm \ref{db:alg:arnoldi}.

\begin{algorithm}[htb]
\begin{algorithmic}[1]
\State{\hspace{\algorithmicindent} 
$[\:{Q}_{1:j-2},\: {w}_{j-1}\:]^T\:[\: {w}_{j-1},\:{Aw}_{j-1}\:]$  \hspace{7mm} \emph{// global reduction}}
\Statex
\State{\hspace{\algorithmicindent} $ {S}_{1:j-2,j} =  {Q}_{1:j-2}^T {Aw}_{j-1}$ \hspace{1.8mm} and \hspace{1mm} $S_{j-1,j} =  {w}_{j-1}^T {Aw}_{j-1}$}
\State{\hspace{\algorithmicindent} $ {C}_{1:j-2,j-1} =  {Q}_{1:j-2}^T {w}_{j-1}$ \hspace{1mm} and \hspace{1mm} $\beta_{j-1} =  {w}_{j-1}^T {w}_{j-1}$}
\Statex
\Statex{\hspace{\algorithmicindent} \emph{// delayed normalization}}
\State{\hspace{\algorithmicindent}$\alpha_{j-1} = \left\{ \:\beta_{j-1} -  {C}^T_{1:j-2,j-1}  {C}_{1:j-2,j-1}\: \right\}^{1/2}$}
\State{\hspace{\algorithmicindent}$T_{j-1,j} = \ddfrac{1}{\alpha_{j-1}^2} \:
\left(\: S_{j-1,j} -  {C}^T_{1:j-2,j-1}\:{S}_{1:j-2,j}\: \right)$}
\State{\hspace{\algorithmicindent}$ {T}_{1:j-2,j} = 
{\ddfrac{1}{\alpha_{j-1}}}\: {S}_{1:j-2,j}$ \label{st_b}}
\Statex
\Statex{\hspace{\algorithmicindent} \emph{// projection}}
\State{\hspace{\algorithmicindent}$ {u}_{j-1} =  {w}_{j-1} -  {Q}_{1:j-2}\:  {C}_{1:j-2,j-1}$}
\State{\hspace{\algorithmicindent}$ {q}_{j-1} = \ddfrac{1}{\alpha_{j-1}}\: {u}_{j-1}$ \label{05}}
\State{\hspace{\algorithmicindent}$ {w}_j = \ddfrac{1}{\alpha_{j-1}}  {Aw}_{j-1} -  {Q}_{1:j-1}\:{T}_{1:j-1,j}$ } 
\Statex
\Statex{\hspace{\algorithmicindent} \emph{// representation $H_{j-1}$}}
\State{\hspace{\algorithmicindent} $ {H}_{1:j-2,j-2} =  {K}_{1:j-2,j-2} +  {C}_{1:j-2,j-1}$}
\State{\hspace{\algorithmicindent} $ {K}_{1:j-1,j-1} =  {T}_{1:j-1,j} - \ddfrac{1}{\alpha_{j-1}}\:H_{1:j-1,1:j-2}\:{C}_{1:j-2,j-1}$} 
\State{\hspace{\algorithmicindent} $H_{j-1,j-2} = \alpha_{j-1}$ }
\end{algorithmic}
\caption{\label{db:alg:arnoldi} Arnoldi-$QR$ (DCGS2)}
\end{algorithm}

\section{Computation and Communication Costs}
The computation and communication costs of the algorithms are listed in  Tables
\ref{db:tbl:iteration} and \ref{db:tbl:full}.  Although theoretically
equivalent, they exhibit different behavior in finite precision arithmetic.
All the schemes, except MGS, are based upon cache-blocked matrix operations.
MGS applies elementary rank--$1$ projection matrices sequentially to a vector
and does not take advantage of the DGEMM matrix-matrix multiplication kernel.
In addition, this algorithm requires one global reduction ({\tt
MPI\_AllReduce}) in the inner-most loop to apply a rank-1 projection matrix.
Thus, $j$ global reductions are required at iteration $j-1$. The implementation
of ICWY-MGS batches together the projections and computes one row of the
strictly lower triangular matrix, \'{S}wirydowicz et
al.~\cite{2020-swirydowicz-nlawa}.
\[
L_{k-1,1:k-2} = \left(\:Q_{1:k-2}^T\:q_{k-1}\:\right)^T.  
\]
The resulting inverse compact $WY$ projector $P$ is given by
\[
P\:a = \left(\: I - Q_{1:j-1}\:T_{1:j-1,1:j-1}\:Q_{1:j-1}^T \:\right)\:a 
\]
where the triangular correction matrix is given by
\[
T_{1:j-1,1:j-1} = (\:I + L_{1:j-1,1:j-1}\:)^{-1}, \quad
T_{1:j-1,1:j-1} \approx (\:Q_{1:j-1}^T\:Q_{1:j-1}\:)^{-1}
\]
The implied triangular solve requires an additional $(j-1)^2$ flops at
iteration $j$$-1$ and thus leads to a slightly higher operation count 
compared to the original MGS orthogonalization scheme, the operation $Q_{1:k-2}^Tq_{k-1}$ increases
ICWY-MGS complexity by $mn^2$ meaning it is 1.5 times more expensive 
($3mn^2$~total) but reduces synchronizations from $j$$-1$ at iteration 
$j$ to 1. This reasoning also follows for the
CWY-MGS algorithm, it is  1.5 times more expensive compared to the sequential 
implementation of MGS.
However, only one global reduction is required per iteration, and hence the 
amount of inter-process communication does not depend upon the number 
of rank--1 projections applied at each iteration.

In the case of the DCGS2 algorithm, the symmetric correction matrix 
$T_{j-1}$ was derived in
Appendix 1 of \cite{2020-swirydowicz-nlawa} and is given by 
\[ 
T_{1:j-1,1:j-1} = I - L_{1:j-1,1:j-1} - L_{1:j-1,1:j-1}^T. 
\]
This form of the projector was employed in the $s$-step and pipe\-lined GMRES
described in Yamazaki et al.~\cite{2020-yamazaki-proceedings-of-siam-pp20}.
When the matrix
$T_{1:j-1,1:j-1}$ is split into $I - L_{1:j-1,1:j-1}$ and $L_{1:j-1,1:j-1}^T$ and 
applied across two iterations of the DCGS2 algorithm, the resulting loss of 
orthogonality is ${\cal O}(\eps)$ in practice.

Block generalizations of the DGCS2 and CGS2 algorithm are presented in Carson
et al.~\cite{2020_arXive_Carson, 2021_arXive_Carson}.  These papers generalize
the Pythagorean trick to block form and derive {\tt BCGS-PIO} and {\tt
BCGS-PIP} algorithms with the more favorable communication patterns described
herein. An analysis of the backward stability of the these block Gram-Schmidt
algorithms is also presented.

\begin{table}[htb]
\begin{center}\footnotesize
\begin{tabular}{|c|c|c|c|}
\hline
orth scheme & flops per iter & synchs & bandwidth \\ \hline\hline
MGS Level 1  & $4(m / p)j$ & $j$ & $j$ \\ \hline    
MGS Level 2 & $6(m / p)j + j^2$ & 1 & $2j$ \\ \hline
CGS & $4(m / p)j$ & 2 & $j$ \\ \hline
CGS2 & $8(m / p)j$ & 3 & $2j$ \\ \hline
CGS2 (lagged norm) & $8(m / p)j$ & 2 & $2j$ \\ \hline
DCGS2-HRT & $8(m / p)j$ & 1 & $2j$ \\ \hline
DCGS2 (QR) & $8(m / p)j$ & 1 & $2j$ \\ \hline
DCGS2 (Arnoldi) & $8(m / p)j + j^2$ & 1 & $2j$ \\ \hline
\end{tabular}
\caption{\label{db:tbl:iteration} 
Cost per iteration for Gram Schmidt algorithms. \\ Where $p$ is the number of processes used.}
\end{center}
\end{table}

\begin{table}[htb]
\begin{center}\footnotesize
\begin{tabular}{|c||c|c|c|} \hline
orth scheme & flops per iter & synchs & bandwidth \\ \hline\hline
MGS Level 1         & $2(m / p)n^2$ & $\frac{1}{2}n^2$ & $\frac{1}{2}n^2$ \\ \hline    
MGS Level 2 & $3(m / p)n^2 + \frac{1}{3}n^3$ & $n$ & $n^2$ \\ \hline
CGS & $2(m / p)n^2$ & $2n$ & $\frac{1}{2}n^2$ \\ \hline
CGS2 & $4(m / p)n^2$ & $3n$ & $n^2$ \\ \hline
CGS2 (lagged norm) & $4(m / p)n^2$ & $2n$ & $n^2$ \\ \hline
DCGS2-HRT & $4(m / p)n^2$ & $n$ & $n^2$ \\ \hline
DCGS2 (QR) & $4(m / p)n^2$ & $n$ & $n^2$ \\ \hline
DCGS2 (Arnoldi) & $4(m / p)n^2 + \frac{1}{3}n^3$ & $n$ & $n^2$ \\ \hline
\end{tabular}
\caption{\label{db:tbl:full} 
Total cost of Gram Schmidt algorithms. \\ Where $p$ is the number of processes used.}
\end{center}
\end{table}

\begin{figure}[htb]
\centering
\includegraphics[width=0.4\textwidth]{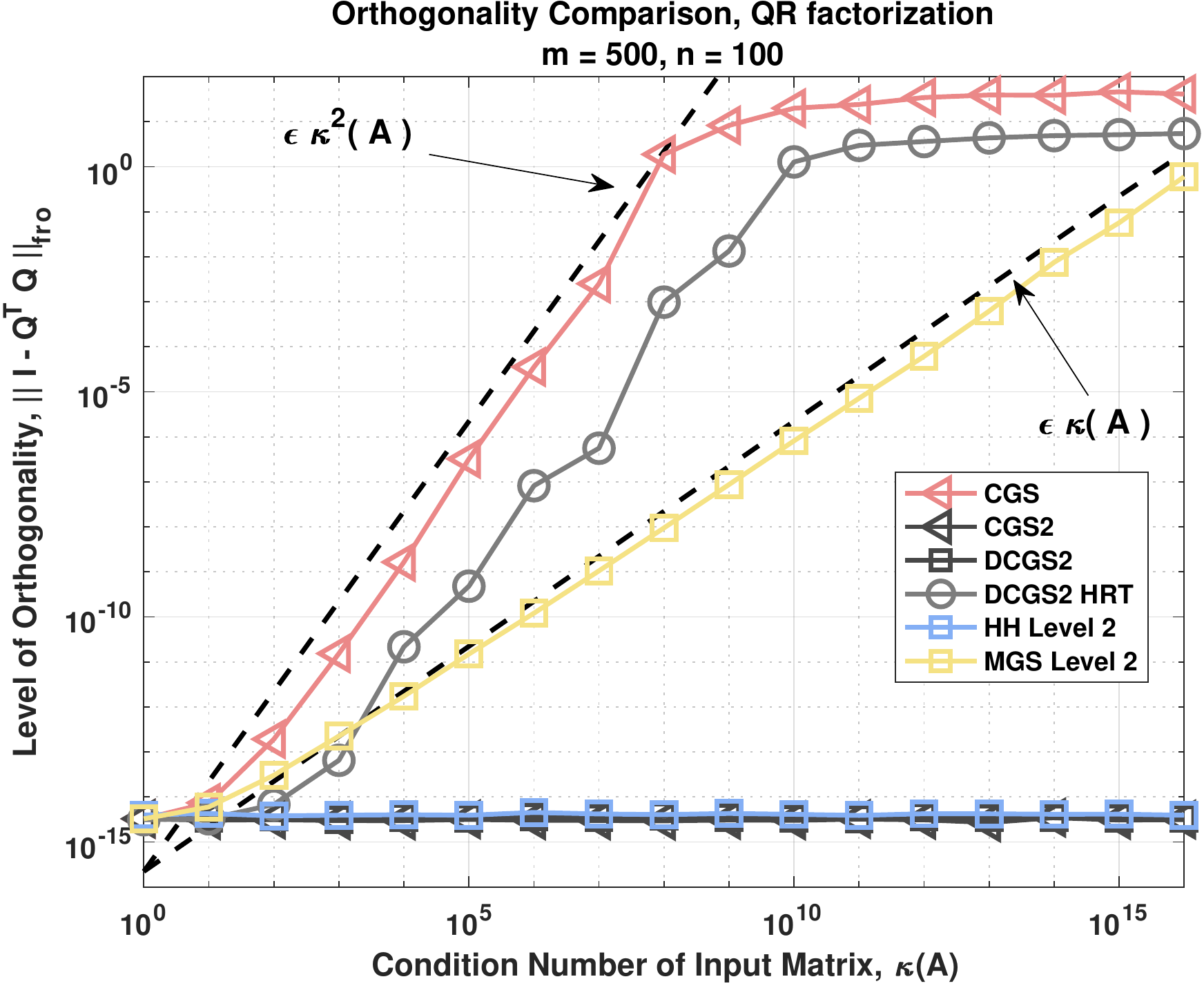}
\caption{\label{fig:orthoerror-cond}Loss of orthogonality with increasing condition number.}	
\end{figure}

\begin{figure}[htb]
\centering
\includegraphics[width=0.4\textwidth]{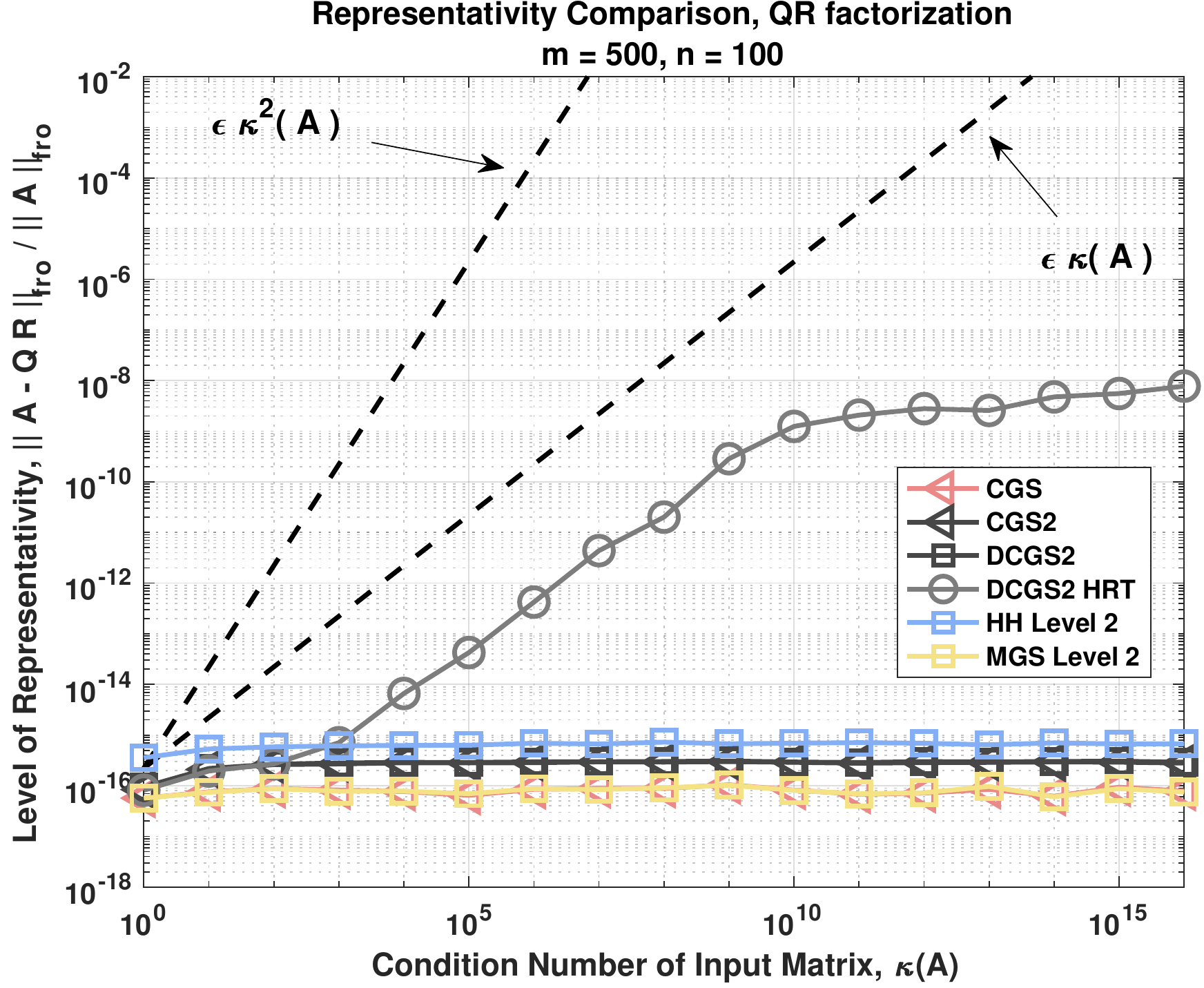}
\caption{\label{fig:represerror-cond} Representation error with increasing condition numbers.}
\end{figure}

It is important to note that there are a variety of ways to implement a blocked
MGS algorithm. The correction matrix $T$ is constructed much like in the 
BLAS Level 2 compact-$WY$ Householder transformation \cite{2020-swirydowicz-nlawa}.
For all results reported herein that employ blocked MGS, these
are based upon CWY-MGS. Except in the context of GMRES. When using the
one-reduce GMRES solver, the underlying Gram-Schmidt algorithm employed is
ICWY-MGS . For this reason Tables \ref{db:tbl:iteration}, \ref{db:tbl:full},
and \ref{db:tbl:loss_of_orth} refer to the blocked MGS implementation as Level
2 versus ICWY-MGS. In addition, MGS Level 1 refers to the sequential implementation for using level 1 BLAS operations.

Both CGS and CGS2 are based upon matrix-vector operations.  CGS applies a
single projection, and then normalizes, requiring two separate steps.  This
projection step consists of two \texttt{DGEMV} kernel calls and one
\texttt{DDOT} for the normalization. CGS suffers from at least an ${\cal
O}(\eps)\kappa^2(A)$ loss of orthogonality.  CGS2 achieves ${\cal O}(\eps)$
through two passes (see Figure \ref{fig:orthoerror-cond}).  The additional
projection within CGS2 accounts for one additional global reduction per
iteration and an additional $4(m/p)j$ operations. 

DCGS2 requires one reduction and employs matrix-matrix multiplies for the
computation in a tall-and-skinny DGEMM.  This leads to the higher sustained
execution rate of DCGS2 (e.g. $2\times$ the GigaFlop/sec). In the context of
Arnoldi, DCGS2 requires an additional $j^2$ flops at the $j$--th iteration.
The additional cost is due to the Arnoldi representation trick described in
Section \ref{db:sec:arnoldi}.  The representation error correction terms
require an additional $n^2$ operations from a matrix-vector product with the
Hessenberg matrix. 

\begin{table}[htb]
\begin{center}\footnotesize
\begin{tabular}{|c||c|c|c|} \hline
orth scheme & LOO & Proven \\ \hline\hline
MGS Level 1 & $\mathcal{O}(\eps) \kappa(A)$ & \cite{Bjork_1967} \\ \hline
MGS Level 2 & $\mathcal{O}(\eps) \kappa(A)$ & Conjectured \\ \hline
CGS & $\mathcal{O}(\eps) \kappa^2(A)$ & \cite{glre:nm:05} \\ \hline
CGS2 & $\mathcal{O}(\eps)$ & \cite{glre:nm:05} \\ \hline
CGS2 (lagged norm) & $\mathcal{O}(\eps)$ & Conjectured \\ \hline
DCGS2-HRT & $\mathcal{O}(\eps)\kappa^2(A)$ & Conjectured \\ \hline
DCGS2 (QR) & $\mathcal{O}(\eps)$ & Conjectured \\ \hline
DCGS2 (Arnoldi) & $\mathcal{O}(\eps)$ & Conjectured \\ \hline
\end{tabular}
\caption{\label{db:tbl:loss_of_orth} Loss of Orthogonality (LOO). }
\end{center}
\end{table}

\section{Numerical Results}\label{num}

In this section, the numerical stability of the Arnoldi algorithm is
investigated for the different 
orthogonalization schemes.
The methodology for the
numerical stability analysis is presented in Section~\ref{sec:manteuffel} along
with the experiments. The same methodology is employed in Section~\ref{ppr}.
Four stability metrics are examined,
\begin{enumerate}
\item representation error
\item loss of orthogonality
\item forward error in the eigenvalue solutions, $<$ threshold
\item dimension of converged invariant subspace, $<$ threshold
\end{enumerate}
The metrics (1) and (2) are sufficient to  analyze the stability of an
orthogonalization scheme. However to give a broader perspective the metrics (3)
and (4) are also examined. Additional metrics that can be considered are:
\begin{enumerate}
\item convergence of GMRES 
\item achievable backward error of GMRES 
\item number of eigenvalue-pairs (Ritz values) with a backward error $<$ threshold 
(see Hernandez et al.~\cite{2007-hernandez-parco})
\end{enumerate}

The convergence of GMRES and the achievable backward error are informative metrics, however, Paige et al.~\cite{2006--simax--paige-rozloznik-strakos} proved that GMRES (with one right-hand side) only needs an
$\mathcal{O}(\eps)\:\kappa(B)$ LOO to converge. Therefore, GMRES is tolerant of ``bad'' 
orthogonalization schemes and is not a stringent enough test.

The number of eigenvalue-pairs with a backward
error less than a threshold should not be used to assess the quality of an
orthogonalization scheme because, for example, a scheme that always returns the same
eigenvalue-pair $n$ times would score the highest possible score ($n$) according to this metric, while performing very poorly in any reasonable metric.

\subsection{Manteuffel Matrix and experimental stability methodology}\label{sec:manteuffel}
The matrix generated by ``central differences" introduced by Manteuffel
\cite{1978_mantueffel_adaptive} is employed in a series of tests designed for
the Krylov-Schur eigenvalue solver based upon DCGS2--Arnoldi.  This matrix is
convenient for computing the forward error solution because the explicit
computation of each eigenvalue is possible, thus the comparison against a
converged eigenvalue is possible.  For $m\times m$ block diagonal matrices $M$
and $N$, where $M$ and $N$ have $k\times k$ sub-blocks such that $m = k^2$.
${M}$ is positive definite and ${N}$ is skew-symmetric. 
An explicit formulation of $M$ and $N$ is given in 
\cite{1978_mantueffel_adaptive}. 
The Manteuffel matrix is expressed as the sum
\begin{equation}
{A} = \frac{1}{h^2}  {M} + \frac{\beta}{2h}  {N}
\end{equation}
where the matrix blocks and $k^2$ eigenvalues are generated by
\begin{equation}
\lambda_{\ell,j} = 2 \left[ 2 - \sqrt{1 - \left(\frac{\beta}{2}\right)^2} \left(\cos\left(\frac{\ell\pi}{L}\right) + 
\cos\left(\frac{j\pi}{L}\right)\right) \right],
\label{eqn:eigs_manteuffel}
\end{equation}
For $ \ell = 1,\dots,k$ and $j = 1,\dots,k$.  As can be seen in
(\ref{eqn:eigs_manteuffel}), $\beta$ is a scalar that governs the spectrum of
the eigenspace.  $L$ is defined by the domain of the differential operator,
$[0,L]\times [0,L]$ and $h=L\,/\,(k+1)$ is the discretization parameter. 
For the experiments within this section, $\beta = 0.5$ and $L=k+1$ so that
$h=1$ ($\beta \leq 2$ implies all eigenvalues are real).  For $k=50$, relevant
numerical metrics are summarized in Table \ref{tab:stab_manteuffel}.  Here, $
{V}$ and $ {W}$ are the left and right eigenvectors.  The Manteuffel matrix is
employed to evaluate the convergence of Krylov-Schur.

\begin{table}[htb!]
\begin{center}
\begin{tabular}{| l | c |} \hline
$\| {A}\|_2$  & $7.99$e$+00$ \\
Cond($ {A}$)  & $3.32$e+$02$ \\
Cond($ {V}$)  & $3.96$e$+11$ \\
Cond($ {W}$)  & $3.74$e$+11$ \\
$\| {A}^T {A} -  {A} {A}^T\|_F / \| {A}\|_F^2$ & $2.81$e$-04$ \\
$\max_{i}$(\, Cond($\;\lambda_{i}\;$) ) & $1.46$e$+10$ \\
$\min_{i}$(\, Cond($\;\lambda_{i}\;$) ) & $2.38$$+02$ \\ 
\hline
\end{tabular}
\caption{\label{tab:stab_manteuffel}Differential operator specs for $k=50$, $m=k^2=2500$. } 
\end{center}
\end{table}

\begin{figure}[htb!]
\centering
\includegraphics[width=0.4\textwidth]{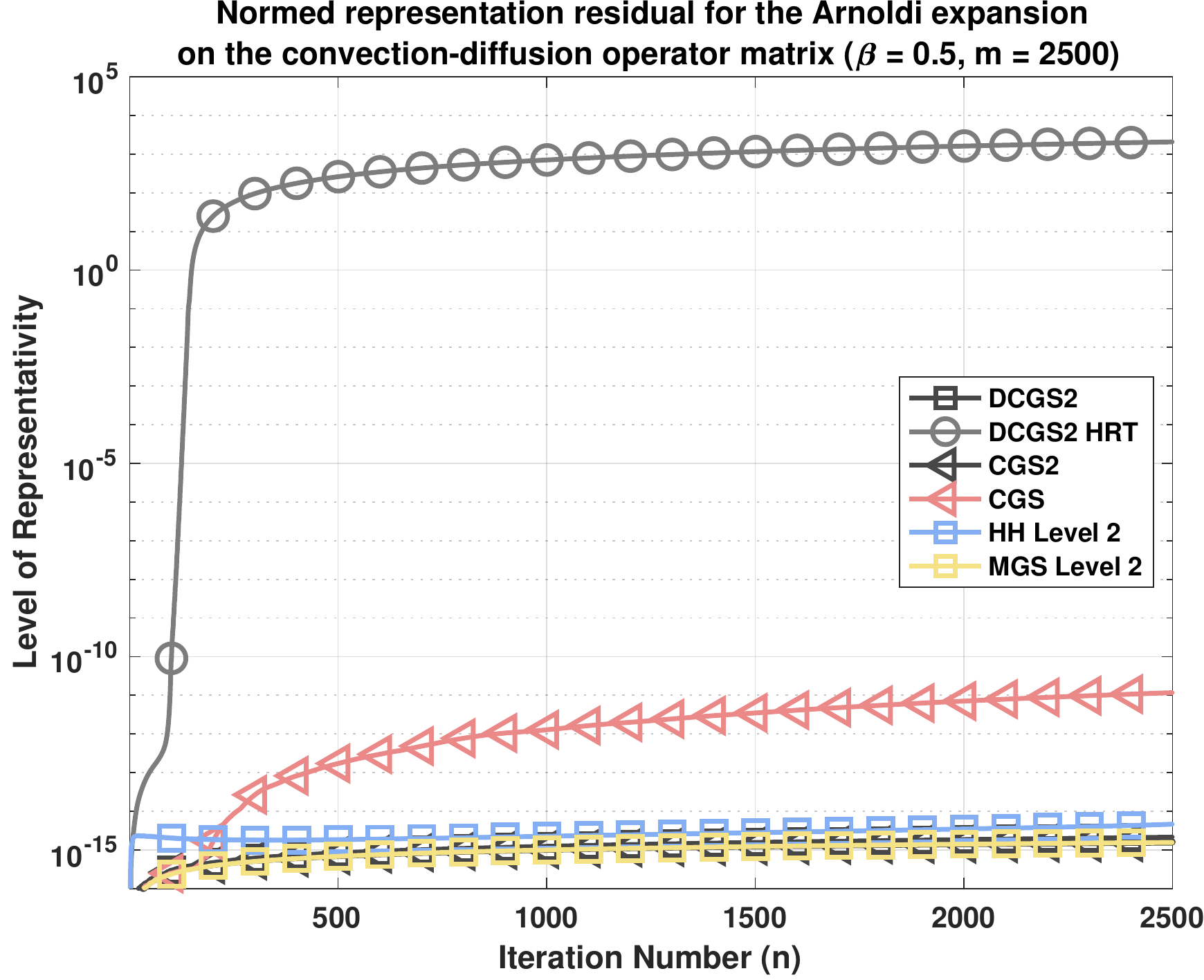}
\caption{\label{fig:stab_manteuffel_repres}
Representation error.}		
\end{figure}

\begin{figure}[htb!]
\centering
\includegraphics[width=0.4\textwidth]{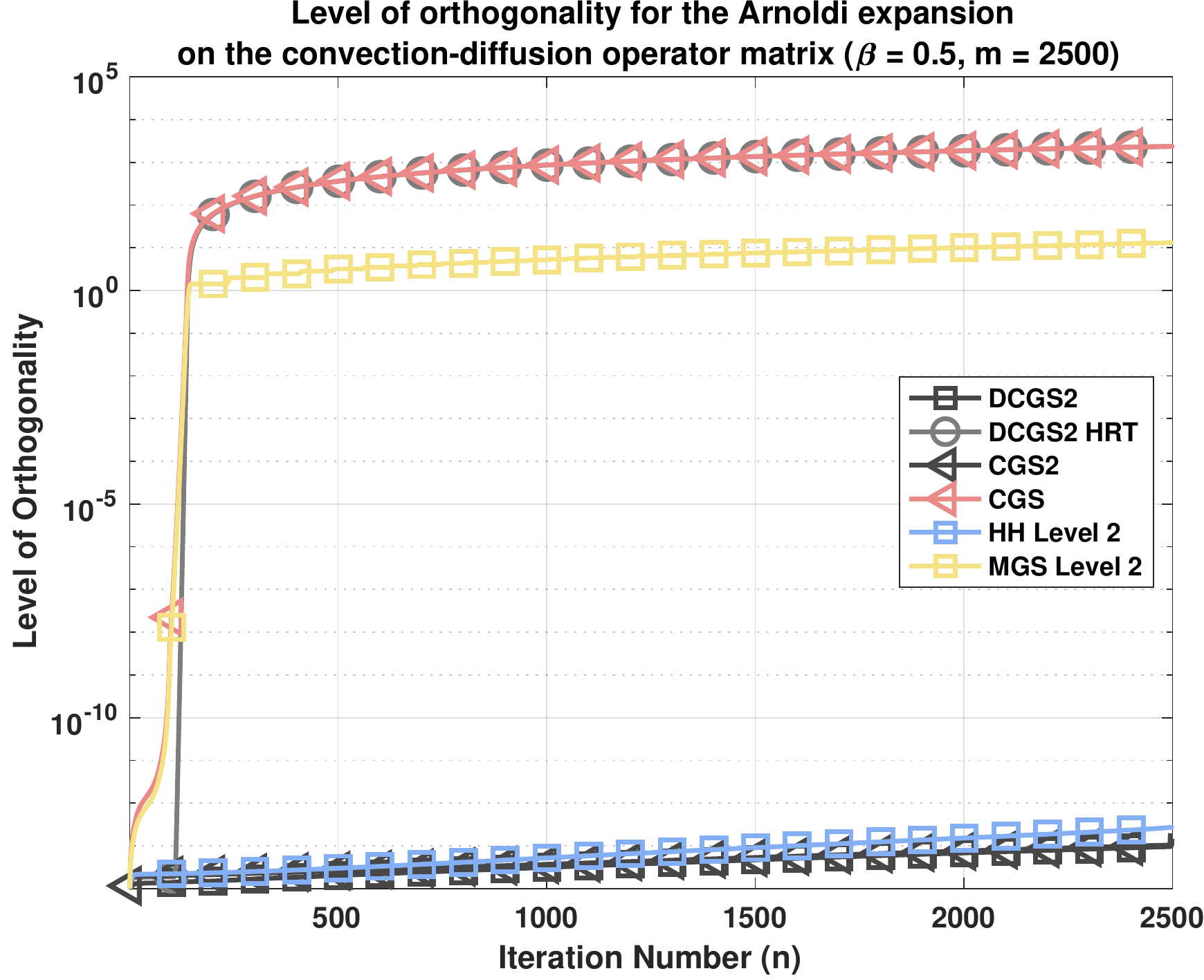}
\caption{\label{fig:stab_manteuffel_orth}
Loss of orthogonality.}		
\end{figure}

The Arnoldi residual and error metrics employed herein are described in
Hern\'{a}ndez et al.~\cite{SLEPc}.  Figure \ref{fig:stab_manteuffel_orth}
displays the loss of orthogonality $\|\: I_{j-1} - {Q}_{1:j-1}^T {Q}_{1:j-1}
\:\|_F$, while Figure \ref{fig:stab_manteuffel_repres} is a plot of the Arnoldi
relative representation error, from (\ref{eq:arnoldi})
\[
RRE(j) = \frac{
\left\|\:{AQ}_{1:j-2} - {Q}_{1:j-1} H_{1:j-1,1:j-2} \:\right\|_F }{ \|  {A} \|_F }.
\]

Each of the algorithms, except for the DCGS2-HRT presented in
\cite{2007-hernandez-parco}, achieve machine precision level relative error.
Figure \ref{fig:stab_manteuffel_orth} displays the loss of orthogonality for
each scheme.  It was noted earlier that CGS exhibits an
$\mathcal{O}(\eps)\kappa^2(A)$ loss of orthogonality. The plot illustrates that
DCGS2-HRT follows CGS while the other algorithms construct orthonormal columns
to the level of machine precision.

\begin{figure}[htb!]
\centering
\includegraphics[width=0.4\textwidth]{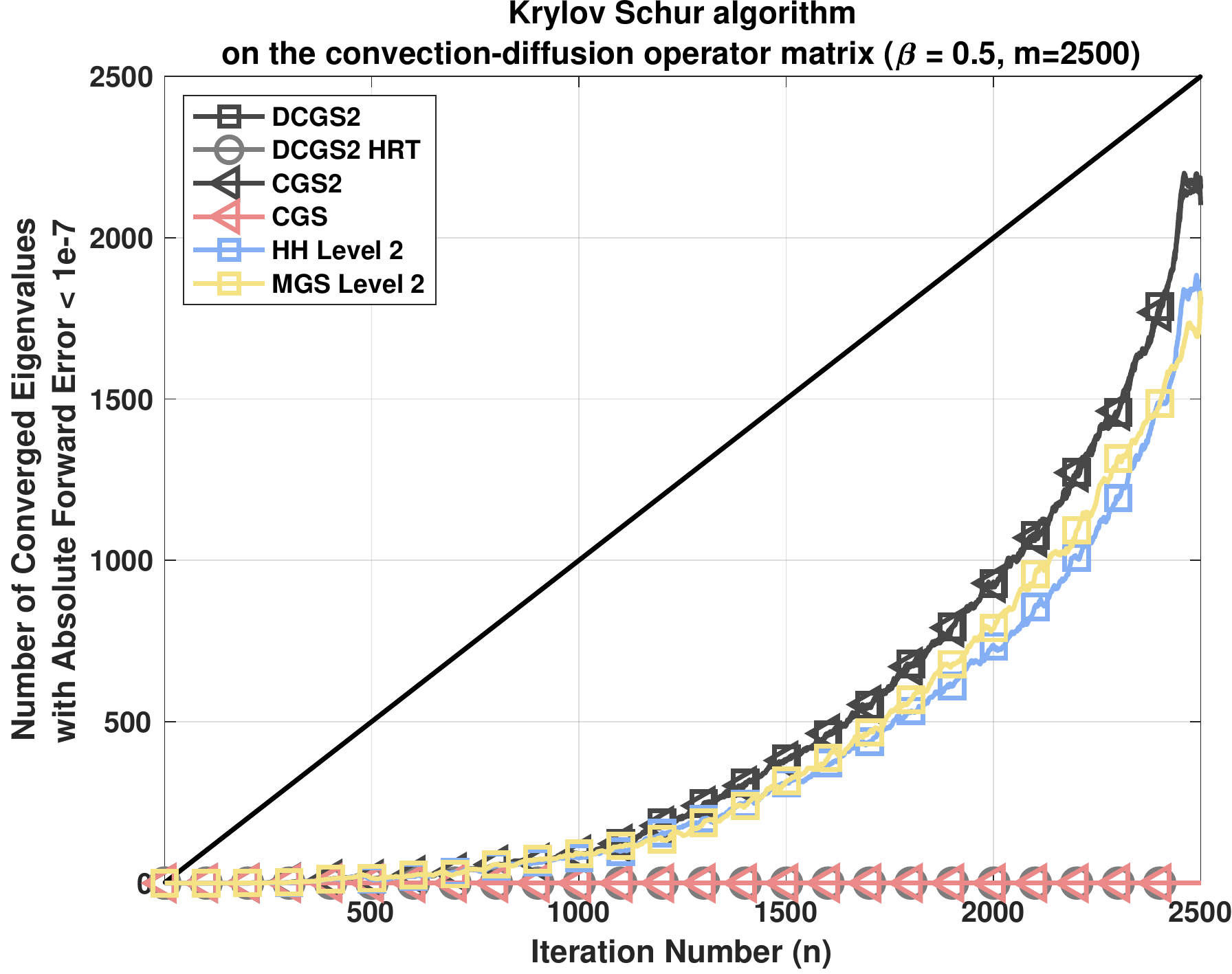}
\caption{\label{fig:stab_manteuffel_eig}Number of eigenvalues computed with $10^{-7}$ absolute error.}	
\end{figure}

\begin{figure}[htb!]
\centering
\includegraphics[width=0.4\textwidth]{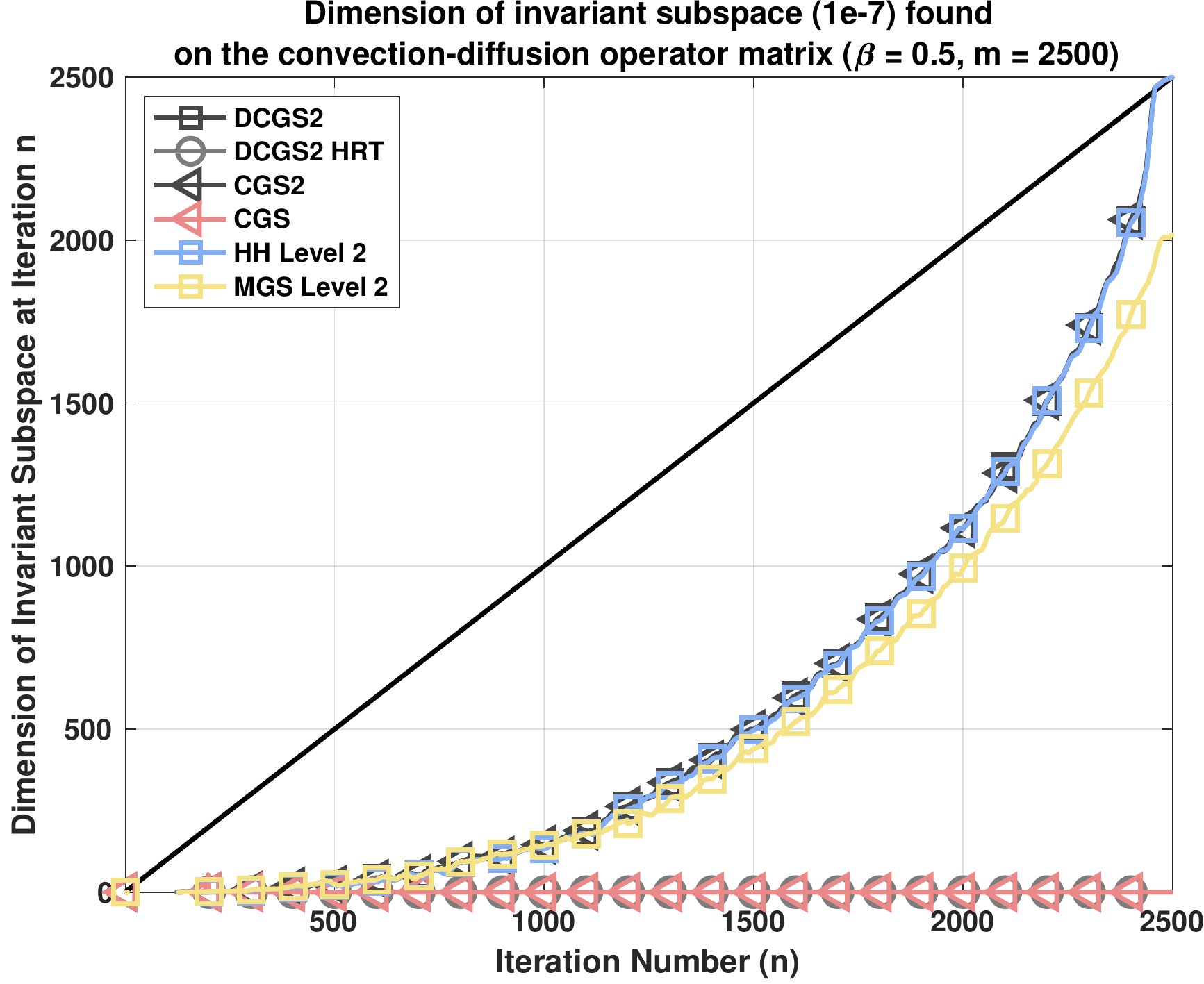}
\caption{\label{fig:stab_manteuffel_iss}Dimension of invariant subspace computed.}
\end{figure}

The results from a Krylov-Schur eigenvalue experiment to evaluate the
convergence properties of the different Arnoldi algorithms are plotted in
Figure \ref{fig:stab_manteuffel_eig}. The solver relies upon the Schur
decomposition of the Hessenberg matrix $H_n$ generated in the Arnoldi
expansion.  To assess the convergence rates, the Arnoldi residual
(\ref{ArnoldiRes}) is compared to the absolute error tolerance.  The
approximate eigenvector (or Ritz vector) associated with the eigenvalue
$\lambda_i$ is defined by $z_i = V_n\:y_i$, where $y_i$ is the corresponding
eigenvector of $H_n$, see  \cite{2007-hernandez-parco}.
\begin{equation}
\label{ArnoldiRes}
\| \: (\: A - \lambda_i\:I \: )\:z_i \:\|_2 = 
|H_{n+1,n}|\: |\: e_n^Ty_i \: | \: < \: \text{tol}
\end{equation}
where ${\rm tol}=1$e$-7$.  If this threshold is satisfied, the iteration is
considered to have found an invariant subspace and the associated diagonal
element in the Schur triangular matrix $T_{l,l}$ is an eigenvalue.  The
representation error and loss of orthogonality can be easily computed.  It is
important to note if these quantities are not close to machine precision, a
converged invariant subspace has not been found.  After the size of the
invariant subspace has been found, the $k^2$ eigenvalues from the formula
in (\ref{eqn:eigs_manteuffel}) are computed and compared with the ``converged"
eigenvalues in the Schur triangular matrix $T$.  In addition, rather than
computing the same eigenvalue twice, the multiplicity is obtained
to determine if the Krylov-Schur algorithm has computed the
same eigenvalue, or unique eigenvalues in the decomposition. The exact
multiplicity of any given eigenvalue was always found.

The plot in Figure \ref{fig:stab_manteuffel_eig} displays the number of
converged eigenvalues at each step $m$ of the Arnoldi algorithm according to
the absolute forward error $|\lambda_{i} - T_{l,l}| <\text{tol}$,
where $\text{tol}=1e-7$. In practice, at each iteration $n$, for each $l$
from $1$ to $n$, the $\lambda_{i}$ are scanned for each $i$ closest to $T_{l,l}$, 
which has not
been found for a previous $l$, and which satisfies the tolerance is selected.
Our code will return an error flag if an iteration returns more eigenvalues
than the expected multiplicity.  The flag was never triggered during our
experiments. At iteration $n$, at most $n$ eigenvalues are found and these
correspond to the solid black line. There are three reasons why the number of
eigenvalues found is not $n$. First, the eigenvalues must be converged. At step
$n$, in exact arithmetic, all eigenvalues would have been found. For step
$n<m$, the number of eigenvalues found is between $0$ and $n$. Second, the
forward error is sensitive to the condition number of the eigenvalues. Some
eigenvalues have condition number of the order $1$e$+10$, (see Table
\ref{tab:stab_manteuffel}), therefore, using $\varepsilon=2.2$e$-16$ accuracy,
our algorithms are not expected to find all eigenvalues at iteration $m=2,500$.
The maximum number of eigenvalues found is about $2,100$ with CGS2 and DCGS2
methods. This condition number problem is present at any restart $n$ and is
intrinsic when using a forward error criteria. Third, the Arnoldi--$QR$
factorization could have errors in fundamental quantities such that the loss of
orthogonality and representation error are large. This may affect the number of
eigenvalues found at iteration $n$. 

Figure \ref{fig:stab_manteuffel_iss} displays, at each restart $n$, the size
of the invariant subspace found.  None of the methods can follow this line, but
as the full Arnoldi expansion of the Manteuffel matrix is approached, any
scheme that maintains orthogonality can continually find new eigenvalues, or new directions to search.
Comparing both plots illustrates that in practice, when eigenvalues are not
known, looking at the size of the invariant subspace can be a good metric. Note
that between the two plots, there is a small gap for the error formula at a
restart of $m=2500$, where this gap is not present in the invariant subspace
plot. The different Arnoldi variants cannot find all of the invariant
subspaces, which is due to the condition number of the eigenvalues.  Comparing
the different $QR$ factorization schemes and the invariant subspace found,
although it loses orthogonality, Arnoldi with MGS can still find new search
directions. Arnoldi--$QR$ based on Householder (HH), CGS2 and DCGS2, can find
a subspace that spans the entire space, but for this matrix MGS still performs
well and generates a subspace size close to 2000. 

\subsection{Matrix Market}

The Arnoldi--$QR$ factorization algorithms are now compared for matrices
gathered from the Suite-Sparse collection maintained by Tim Davis at Texas A\&M
University \cite{suite_sparse_collection}. A total of 635 matrices were chosen
by the following criteria: (1) number of nonzeros $< \;500,000$, (2) the matrix
is  REAL, (3) the matrix is UNSYMMETRIC and (4) the number of columns and rows
$> \;100$. The Krylov basis is computed for each of the 635 matrices in the 
collection. The representation error and loss of orthogonality are computed
for every 5 columns until 75 (making sure the dimension of any matrix is not
exceeded). Meaning the restart in an Arnoldi expansion varies from $n=5$ to
$n=75$ in increments of $5$. 

\begin{figure}[htb!] 
\centering
\includegraphics[width=0.4\textwidth]{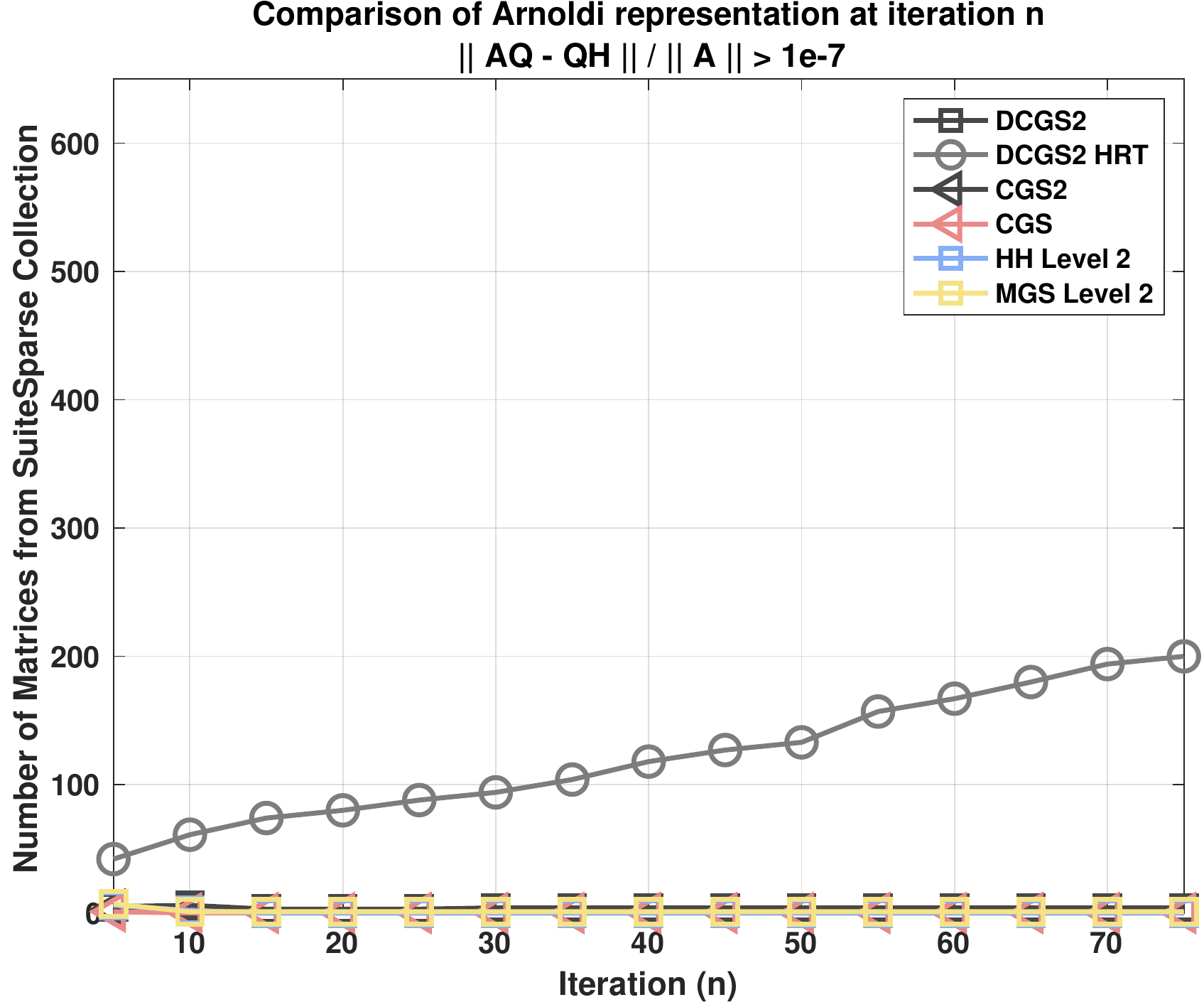}
\caption{\label{fig:635_suitesparse_repres}Representation error for Suite-Sparse
matrices.}	
\end{figure}

\begin{figure}[htb!]
\centering
\includegraphics[width=0.4\textwidth]{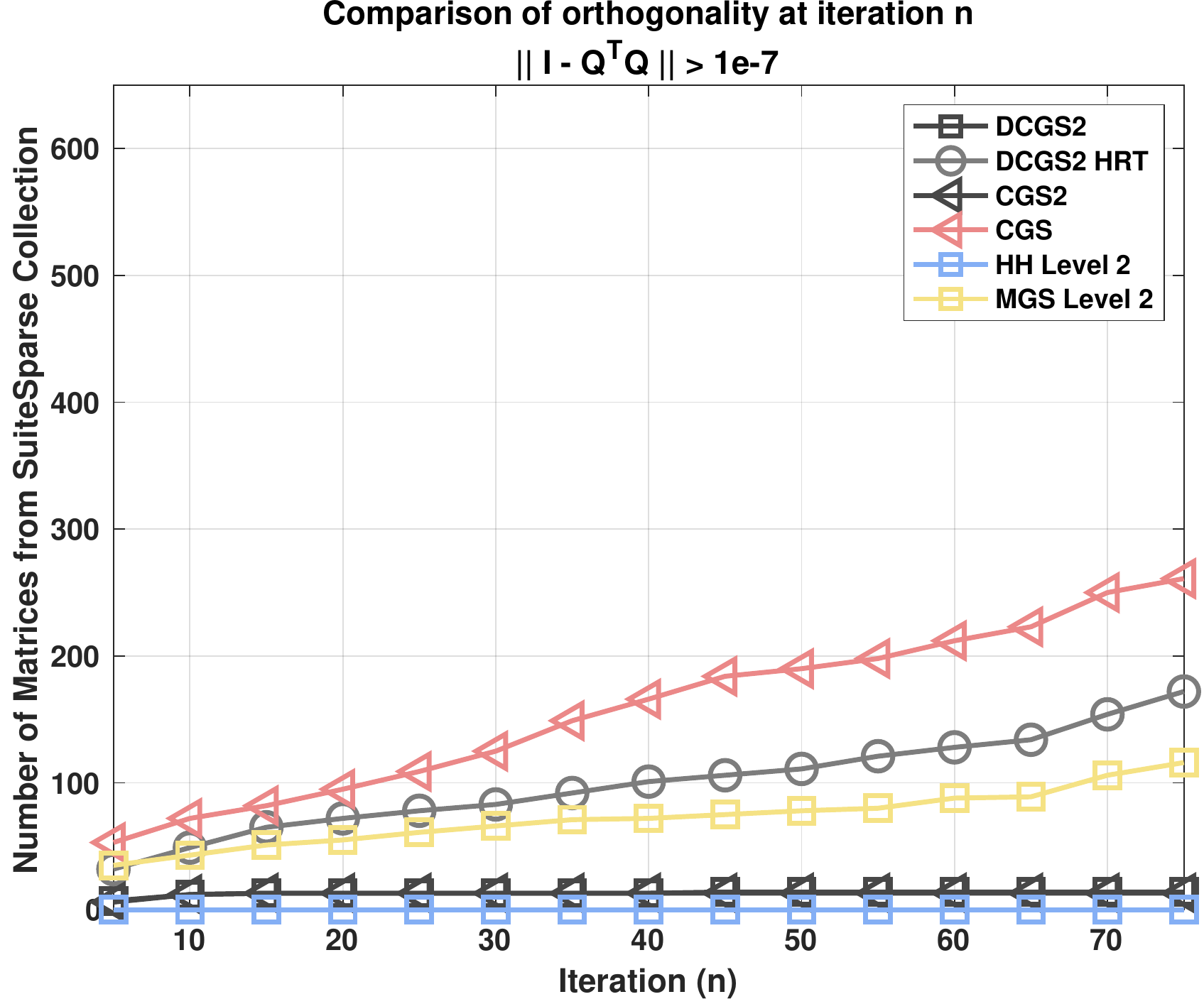}
\caption{\label{fig:635_suitesparse_orth}Loss of Orthogonality for Suite-Sparse matrices.}	
\end{figure}

Figures \ref{fig:635_suitesparse_repres} and \ref{fig:635_suitesparse_orth}
display these metrics for each of the schemes. At each iteration the tolerance
is set to $1$e$-7$. If the representation error or loss of orthogonality is
above this threshold the matrix is flagged.  The $y$-axis represents the total
number of matrices above the given threshold and the $x$-axis indicates the
Krylov subspace dimension (restart $m$) employed by the Arnoldi expansion. 

Figure \ref{fig:635_suitesparse_orth} clearly indicates that Krylov vectors
generated using CGS and MGS lose orthogonality at different rates.  It is
observed that the DCGS2-HRT curve falls between these.  For the Manteuffel
matrix, DCGS2-HRT  appears to perform more like CGS and lies somewhere in
between.  It is important to note, in Figure \ref{fig:635_suitesparse_repres},
that DCGS2-HRT does not maintain a low representation error for the Arnoldi
expansion. This is also apparent in Figure \ref{fig:represerror-cond}.

With a restart of $n=75$, these metrics are plotted in Table
\ref{tab:635_matrixmarket}.  The additional metric displayed is the size of the
invariant subspace found, described in Section \ref{sec:manteuffel}.

\begin{table}[htb!]
\begin{center}
\begin{tabular}{| c | c | c | c |} 
\hline
Orth.        & Repres    & LOO       & Invariant  \\
scheme       & $<$ $1$e$-7$  & $<$ $1$e$-7$  & Subspace   \\
\hline                                            
DCGS2	           & 631             & 621           & \textbf{9844}      \\
DCGS2 HRT	   & 435             & 463           & 7168               \\
CGS2	           & \textbf{635}    & 622           & 9677               \\
CGS		   & \textbf{635}    & 374           & 8370               \\
HH Level 2	           & \textbf{635}    & \textbf{635}  & 9783               \\
MGS Level 2	           &  634  &  519 &  9580         \\
\hline
\end{tabular}
\caption{\label{tab:635_matrixmarket} $n$ = 75; $\text{tol} = 1$e$-7$; Suite-Sparse matrices} 
\end{center}
\end{table}

\section{Parallel Performance Results}\label{ppr}
Parallel performance results are now presented for the Summit Supercomputer at
Oak Ridge National Laboratory.  Each node of Summit consists of two 22-core IBM
Power 9 sockets and six NVIDIA Volta 100 GPUs.  CGS2 and DCGS2 were implemented
and tested using the Trilinos-Belos iterative solver framework \cite{Belos,
heroux_trilinos_2005}.  Therefore, although NVIDIA V100 results are
presented here, the implementation is portable to different hybrid node
architectures with a single code base. 

To summarize, DCGS2 achieves faster compute times than CGS2 for two reasons.
First, the former employs either matrix-vector or matrix-matrix kernels, which
provide greater potential for data reuse.  For tall-and-skinny matrices,
employed by DGEMV and DGEMM, compute time is often limited by data movement,
and matrix-matrix type kernels often achieve faster execution rates.
Therefore, DCGS2 is faster than CGS2, even on a single GPU.  Second, on
multiple GPUs, the low-synch algorithm decreases the number of
global-reductions.  Therefore, a greater speedup is achieved on a large number
of GPUs.  In this section, the execution rates on a single and multiple GPUs
are compared.

\subsection{Single GPU Performance}\label{sec:perf-single}

Figure \ref{tab:perf_comparison_1gpu} provides the execution rates in 
GigaFlops/sec of the main computational kernels on a single GPU, with an
increasing number of rows or columns, as reported in 
\ref{tab:cgs2_dcgs2_perf_comparison_row_exp_1gpu} and columns~\ref{tab:cgs2_dcgs2_perf_comparison_col_exp_1gpu}
respectively.  Within the plot, 
\begin{itemize}
\item {\tt MvTransMv} computes the dot-products, e.g., 
      \texttt{DGEMV} to compute $S_{1:j-1,j} = Q_{1:j-1}^T a_j$ in CGS2 or 
      \texttt{DGEMM} to compute $[\:Q_{1:j-2},\: w_{j-1}]^T [\:w_{j-1},\: a_j\:]$ in DCGS2.
\item {\tt MvTimesMatAddMv} updates the vectors by applying the projection, e.g.,
      \texttt{DGEMV} to compute $w_j = q_j - Q_{1:j-1}\:S_{1:j-1,j}$ in DCGS2 or
      \texttt{DGEMM} to compute 
\[ 
\left[\:u_{j-1},\: w_{j}\:\right] = 
\left[\:w_{j-1},\:a_j\:\right] - Q_{1:j-2}\: \left[\:C_{1:j-2,j-1},\: S_{1:j-1,j}\:\right] 
\]
\item {\tt MvDot} computes \texttt{DDOT} product of two vectors,
       and is used to compute the normalization factor $\alpha_{j-1} = \|u_j\|_2$.
\end{itemize}
\begin{figure*}[htb]
\begin{center}\footnotesize
\begin{subfigure}[b]{0.45\textwidth}
        \begin{tabular}{| l || c | c | c | c | c |}
        \hline
        Orthog  & \multicolumn{5}{c|}{Number of rows, $n$, in millions}\\
        Scheme  & $1$ & $2$ & $4$ & $8$ & $16$  \\
        \hline
        \hline
        DCGS2           &       &       &       &       & \\
        \; MVTimes GF/s & 300.3 & 319.1 & 302.3 & 320.3 & 331.4 \\
        \; MVTrans GF/s & 201.4 & 211.8 & 191.4 & 150.2 & 126.8 \\
        \; Total GF/s    & 215.1 & 232.3 & 218.2 & 193.3 & 174.8 \\
        \hline
        CGS2       &       &       &       &       & \\
        \; MVTimes GF/s & 132.0 & 136.4 & 153.4 & 163.8 & 169.5 \\
        \; MVTrans GF/s & 128.5 & 141.0 & 146.4 & 135.3 & 122.2 \\
        \; Total GF/s    & 126.4 & 135.4 & 146.7 & 145.5 & 139.6 \\
        \hline
        \end{tabular}
        \caption{\label{tab:cgs2_dcgs2_perf_comparison_row_exp_1gpu} 
                 Fixed number of columns $n=50$.}
\end{subfigure}
\begin{subfigure}[b]{0.45\textwidth}
	\begin{tabular}{| l || c | c | c | c | c |} 
	\hline
        Orthog   & \multicolumn{5}{c|}{Number of columns, $n$}\\
	Scheme  & $100$ & $120$ & $140$ & $160$ & $180$  \\
	\hline                                            
	\hline                                            
        DCGS2        &       &        &        &        & \\
	\; MVTimes GF/s & 353.7  & 362.8  & 369.8  & 375.6  & 379.7 \\ 	
	\; MVTrans GF/s & 169.2  & 164.3  & 160.3  & 158.7  & 155.6 \\
	\; Total GF/s    & 221.2  & 219.8  & 218.4  & 218.5  & 216.6 \\
        \hline
        CGS2       &       &       &       &       & \\
	\; MVTimes GF/s & 182.1  & 186.9  & 190.7 & 193.1  & 195.7  \\ 	
	\; MVTrans GF/s & 152.0  & 153.6  & 153.4  & 153.3  & 153.2  \\ 	
	\; Total GF/s    & 163.8  & 167.1  & 168.6  & 169.7  & 170.8  \\
	\hline
	\end{tabular}
	\caption{\label{tab:cgs2_dcgs2_perf_comparison_col_exp_1gpu} 
	         Fixed number of rows $m=5$e$+6$.} 
\end{subfigure}
\caption{\label{tab:perf_comparison_1gpu}Execution rate (GigaFlops/sec) of BLAS kernels  (1 node, 1 GPU).}
\end{center}
\end{figure*}

\begin{figure}[htb]
        \centering
        \includegraphics[width=0.4\textwidth]{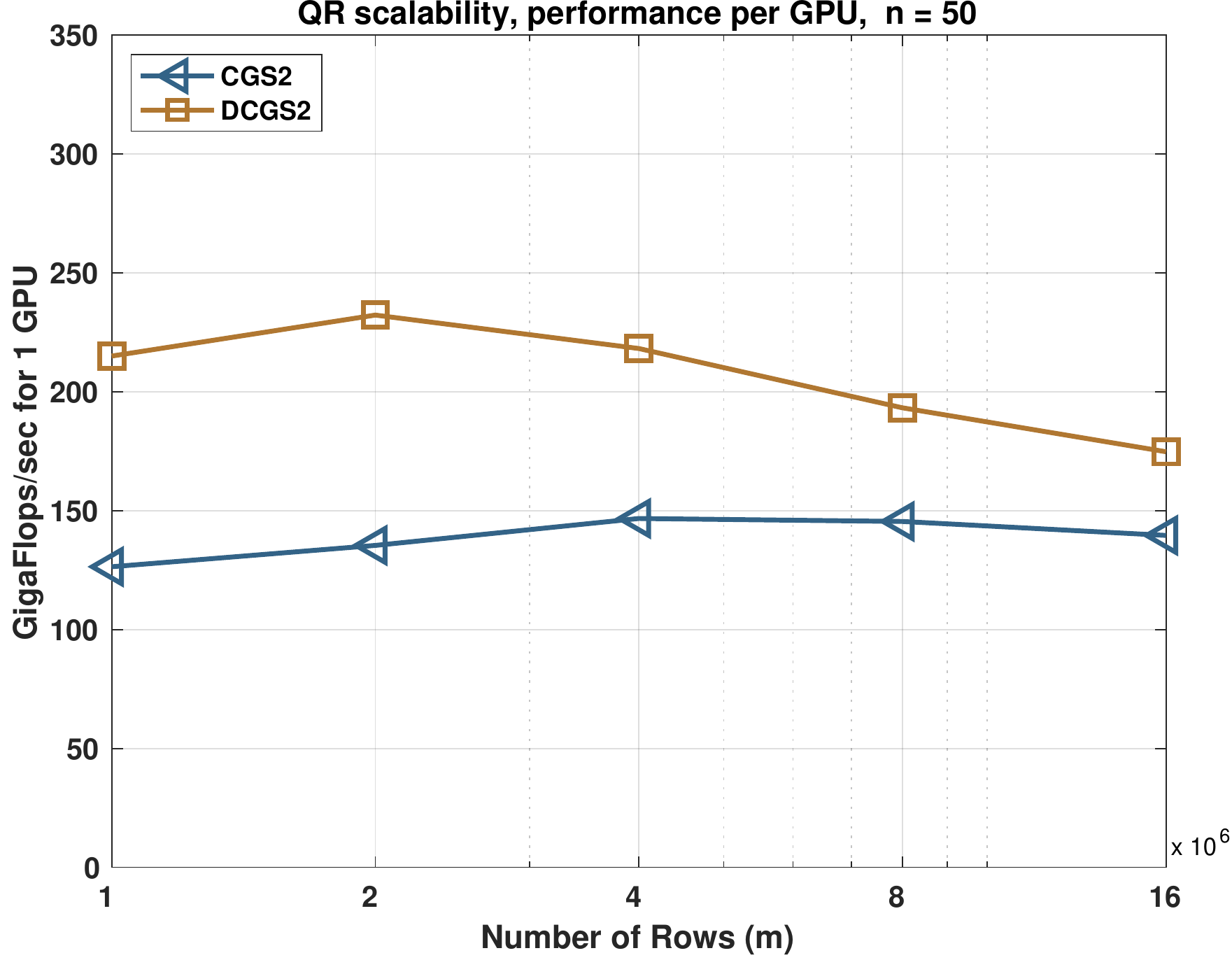}
        \caption{\label{fig:row_exp_compare_1gpu_perf}
        Execution rate (GigaFLops/sec, 1 GPU). $n=50$ columns.}
\end{figure}

\begin{figure}[htb]
        \centering
        \includegraphics[width=0.4\textwidth]{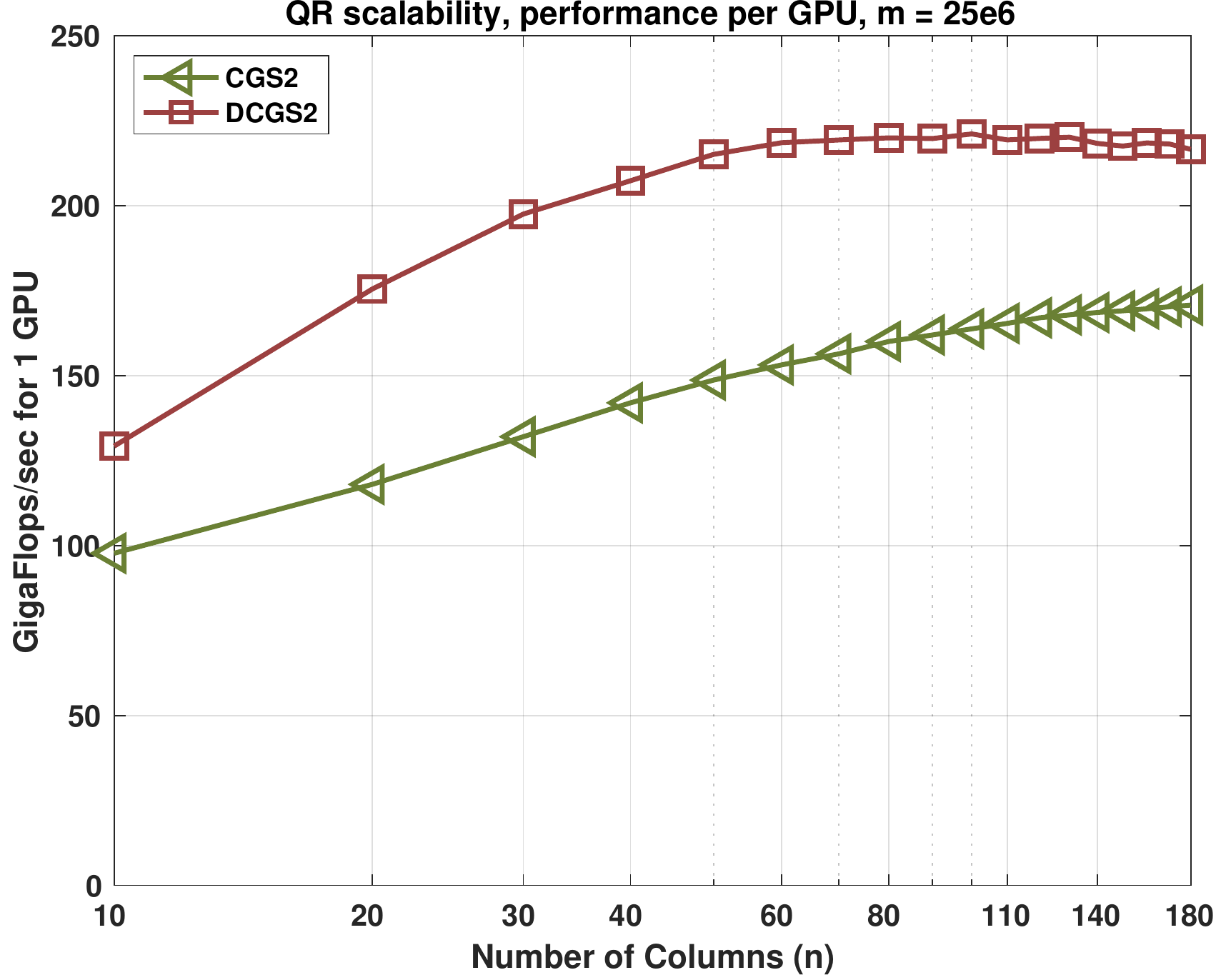}
        \caption{\label{fig:col_exp_compare_1gpu_perf} 
        Execution rate (GigaFlops/sec, 1 GPU). $m=25$e$+6$ rows.} 
\end{figure}

Memory bandwidth utilization is a predictor of performance.  For example, at
the $j$--th iteration of CGS2, {\tt MvTransMV} reads the $m\times (j-1)$ matrix
$Q_{1:j-1}$ and the input vector $a_j$ of length $m$, then writes the result
back to the output vector $S_{1:j-1,j}$, while performing $(2m-1)\times (j-1)$
flops with two flops per read, assuming the vectors remain in caches, or two
flops per two reads and one write with a read-write of vector elements for each
multiply-add.  Thus, on the V100 with a memory bandwidth of 800 GB/s, 200
GigaFlops/sec is expected from this kernel in double precision.  Figures
\ref{fig:row_exp_compare_1gpu_perf} and 
\ref{fig:col_exp_compare_1gpu_perf} display kernel
compute times on one Summit node using one GPU as the number of rows or columns
is varied.

Note that DCGS2 combines two {\tt MvTransMv} calls with a single input and
output vector into a call to {\tt MvTimesMatAddMv} with two input and output
vectors. This can double the potential peak speed (i.e. the $m\times (j-1)$ matrix
is read only once to perform four flops per read).  Figure
\ref{tab:perf_comparison_1gpu}
indicates that for a large number of rows or columns,
that DCGS2 increases the execution rate by up to $1.7$ and $1.4\times$,
respectively..

\subsection{Strong-Scaling Performance}\label{sec:perf-parallel}

\begin{figure}[htb]
        \centering
        \includegraphics[width=0.4\textwidth]{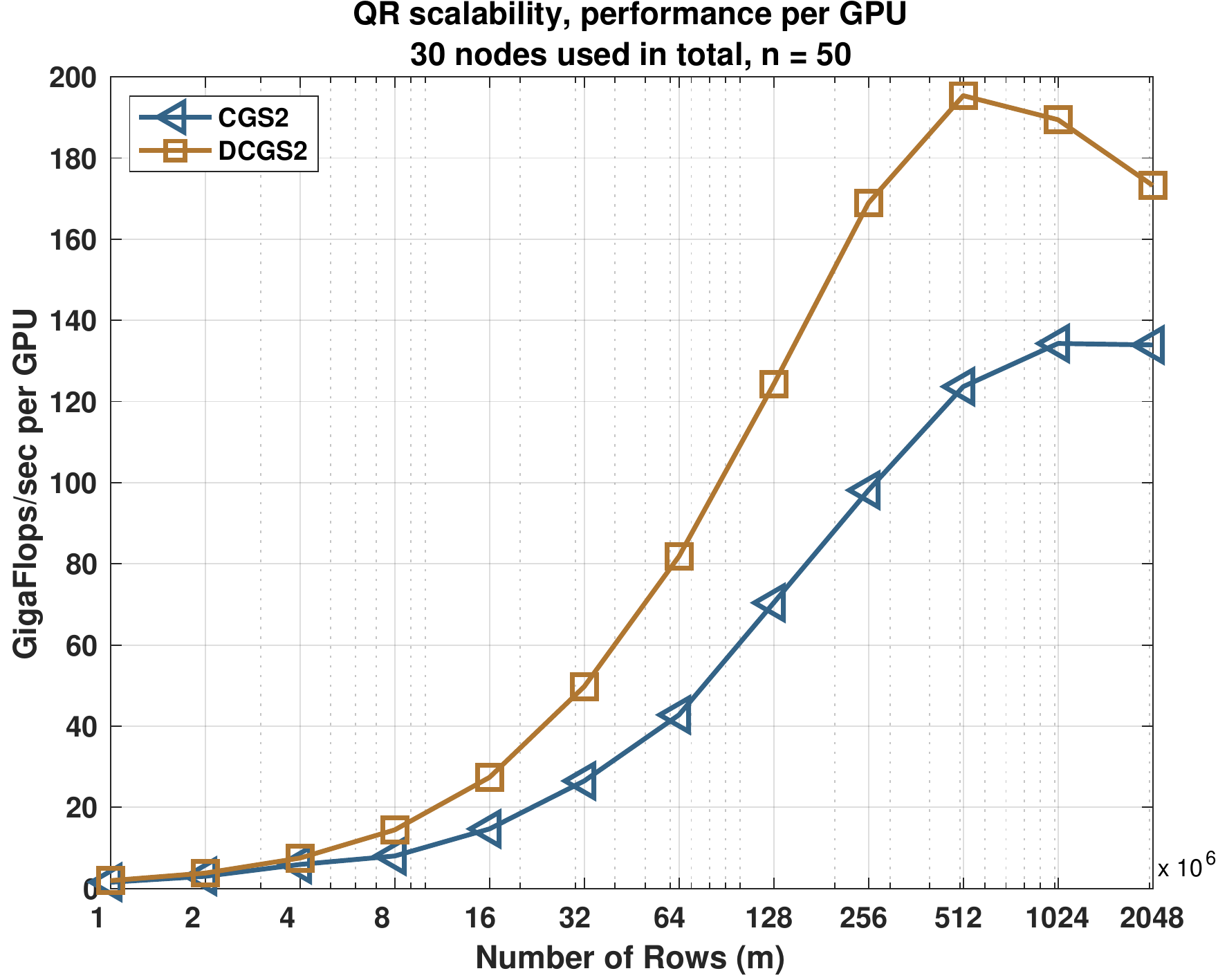}
        \caption{\label{fig:row_exp_compare}
        Execution rate per node (30 nodes, 6 GPUs per node). }
\end{figure}

\begin{figure}[htb]
        \centering
        \includegraphics[width=0.4\textwidth]{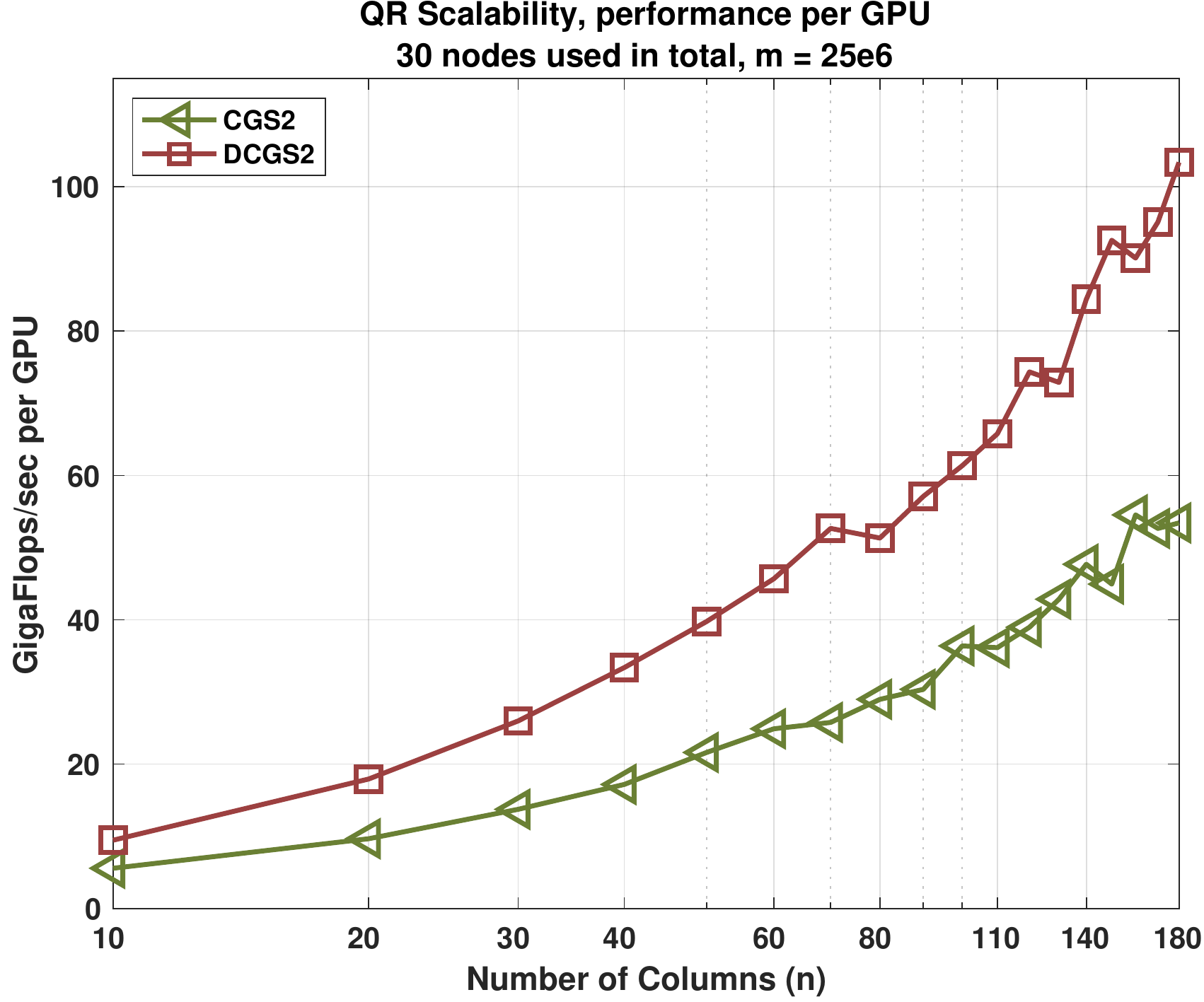}
        \caption{\label{fig:col_exp_compare}
        Execution rate per node (30 nodes, 6 GPUs per node).}
\end{figure}

The speedups obtained by DCGS2 for the two main kernels are presented
in Tables~\ref{tab:cgs2_dcgs2_perf_comparison_communication},
\ref{tab:cgs2_dcgs2_perf_comparison_sequential} and
\ref{tab:cgs2_dcgs2_perf_comparison_total}, while
Figures~\ref{fig:row_exp_compare} and \ref{fig:col_exp_compare} and 
displays the GigaFlops/sec execution rate achieved by the kernels
on 30 nodes, using 6 GPUs per node. Figures
\ref{fig:perf-time-belos-750} and \ref{fig:perf-time-belos-1000} 
represent a strong-scaling study and display the time to solution 
while varying the number of GPUs for a fixed matrix size. One 
run or trial for each node count was employed for the algorithms 
in order to collect the performance data on Summit.

\begin{itemize}
\item
Table \ref{tab:cgs2_dcgs2_perf_comparison_communication} displays the speedup 
(ratio of DCGS2 to CGS2 compute time) for {\tt MvTransMv}.
Because DCGS2 employs fewer global reductions, as the number of GPUs increases,
the speedup obtained by {\tt MvTransMv} increases, reaching up to $2.20\times$ 
faster times on 192 GPUs.
\item
Table \ref{tab:cgs2_dcgs2_perf_comparison_sequential} displays the 
speedup for the {\tt MvTimesMatAddMv} kernel. 
DCGS2 merges two {\tt MvTransMv}
calls into one {\tt MvTimesMatAddMv} and achieves $2\times$ 
speedup on a single GPU.  
With more GPUs, the number of local rows and speedup 
decrease.  However, the compute time is dominated by the 
{\tt MvTimesMatAddMv} kernel.
\end{itemize}
Table \ref{tab:cgs2_dcgs2_perf_comparison_total} displays the speedup
obtained by DCGS2, when varying the number of rows. By combining
matrix-vector products with global reductions, the speedup obtained by DCGS2 
in some instances was significant, up to $3.6\times$ faster.


Figures \ref{fig:perf-time-belos-750} and \ref{fig:perf-time-belos-1000}
display the time to solution for the GMRES linear solvers.
The latter achieves improved strong scaling due to the merged 
{\tt MvTimesMatAddMv} kernel.

\begin{figure*}[hbt]
\begin{center}\footnotesize
\begin{subfigure}[b]{0.45\textwidth}
        \begin{tabular}{| l || c | c | c | c | c |}
        \hline
        Orthog  & \multicolumn{5}{c|}{Number of rows, $n$, in millions}\\
        Scheme  & $128$ & $256$ & $512$ & $1024$ & $2048$  \\
        \hline
        \hline
        DCGS2           &       &       &       &       & \\
        \; MVTimes GF/s & 194.2 & 248.9 & 285.9 & 312.5 & 324.2 \\
        \; MVTrans GF/s & 129.4 & 170.4 & 179.9 & 154.1 & 128.6 \\
        \; Total GF/s    & 124.3 & 168.9 & 195.4 & 189.4 & 173.2 \\
        \hline
        CGS2       &       &       &       &       & \\
        \; MVTimes GF/s & 99.6  & 126.6 & 145.1 & 158.9 & 167.1 \\
        \; MVTrans GF/s & 73.9  & 100.8 & 124.9 & 128.4 & 119.3 \\
        \; Total GF/s   & 42.8  & 98.2  & 123.7 & 134.3 & 133.9 \\
        \hline
        \end{tabular}
        \caption{\label{tab:cgs2_dcgs2_perf_comparison_row_exp} 
                 Fixed number of columns $n=50$. GigaFlops/sec }
\end{subfigure}
\begin{subfigure}[b]{0.45\textwidth}
	\begin{tabular}{| l || c | c | c | c | c |} 
	\hline
    Orthog  & \multicolumn{5}{c|}{Number of columns, $n$}\\
	Scheme  & $100$ & $120$ & $140$ & $160$ & $180$  \\
	\hline                                            
	\hline                                            
        DCGS2        &       &       &       &       & \\
	\; MVTimes GF/s & 121.6 & 135.5 & 151.2 & 159.7 & 168.9 \\ 	
	\; MVTrans GF/s & 61.1  & 76.9  & 89.1  & 89.1  & 107.6 \\
	\; Total GF/s   & 61.4  & 74.4  & 84.4  & 90.1  & 103.4 \\
        \hline
        CGS2       &       &       &       &       & \\
	\; MVTimes GF/s & 122.2 & 95.2  & 101.3 & 99.3  & 90.4  \\ 	
	\; MVTrans GF/s & 38.1  & 39.4  & 49.6  & 56.5  & 52.6  \\ 	
	\; Total GF/s    & 36.4  & 38.9  & 47.7  & 54.5  & 53.4  \\
	\hline
	\end{tabular}
	\caption{\label{tab:cgs2_dcgs2_perf_comparison_col_exp} 
	         Fixed number of rows $n=25$e$+6$. GigaFlops/sec} 
\end{subfigure}
\caption{\label{tab:perf_comparison}Execution rate (GigaFlops/sec) of BLAS kernels 
(30 nodes, 6 GPUs per node ).}
\end{center}
\end{figure*}

\begin{table}[H]
\begin{center}
	\begin{tabular}{| c || c | c | c | c | c |} 
	\hline
                  & \multicolumn{5}{c|}{Number of rows, $n$}\\
	\# GPUs   & $1$e$+6$ & $5$e$+6$ & $10$e$+6$ & $25$e$+6$ & $50$e$+6$  \\
	\hline                                            
	6	  & 1.7        & 1.6        & 1.6         & 1.3         & 1.1 \\
	12	  & 1.9        & 1.8        & 1.7         & 1.5         & 1.3 \\
	24	  & 2.1        & 1.7        & 1.7         & 1.6         & 1.5 \\
	48	  & 2.1        & 1.9        & 1.9         & 1.8         & 1.6 \\
	96    & 2.0        & 2.1        & 4.0         & 1.8         & 1.8 \\
	192   & 2.2        & 2.1        & 2.3         & 2.1         & 2.2 \\
	\hline
	\end{tabular}
	\caption{\label{tab:cgs2_dcgs2_perf_comparison_communication} 
	Speedup of DCGS2 over CGS2 for MvTransMv and MvDot.} 
\end{center}
\end{table}

\begin{table}[htb]
\begin{center}
	\begin{tabular}{| c || c | c | c | c | c |} 
	\hline
                  & \multicolumn{5}{c|}{Number of rows, $n$}\\
	\# GPUs   & $1$e$+6$ & $5$e$+6$ & $10$e$+6$ & $25$e$+6$ & $50$e$+6$  \\
	\hline                                            
	6	  & 0.8       & 2.0        & 2.0         & 1.9         & 2.0 \\
	12	  & 1.9        & 2.0        & 1.9         & 2.0         & 1.9 \\
	24 	  & 1.3        & 0.9        & 1.2         & 2.0         & 2.0 \\
	48	  & 1.1        & 0.9        & 0.9         & 1.9         & 2.0 \\
	96    & 1.1        & 0.8        & 5.2         & 1.1         & 1.3 \\
	192  & 0.9        & 0.7        & 0.7         & 0.9         & 1.3 \\
	\hline
	\end{tabular}
	\caption{\label{tab:cgs2_dcgs2_perf_comparison_sequential} Speedup of DCGS2 over CGS2 for MvTimesMatAddMv.} 
\end{center}
\end{table}

\begin{table}[htb]
\begin{center}
	\begin{tabular}{| c | c | c | c | c | c |} 
	\hline
                  & \multicolumn{5}{c|}{Number of rows, $n$}\\
	\# GPUs   & $1$e$+6$ & $5$e$+6$ & $10$e$+6$ & $25$e$+6$ & $50$e$+6$  \\
	\hline                                            
	6	  & 1.1        & 1.6        & 1.7         & 1.5         & 1.3 \\
	12	  & 1.1        & 1.8        & 1.7         & 1.6         & 1.5 \\
	24	  & 1.1        & 1.5        & 1.8         & 1.7         & 1.6 \\
	48	  & 1.2        & 1.1        & 1.8         & 1.8         & 1.7 \\
	96    & 1.2        & 1.3        & {\bf 4.0}   & 1.8         & 1.8 \\
	192   & 1.4        & 1.3        & 1.7         & 1.9         & {\bf 2.1} \\
	\hline
	\end{tabular}
	\caption{\label{tab:cgs2_dcgs2_perf_comparison_total} Overall speedup 
	of DCGS2 over CGS2.} 
\end{center}
\end{table}

Figure \ref{tab:perf_comparison} displays the GigaFlops/sec execution rates
obtained by CGS2 and DCGS2, along with the BLAS kernels for a fixed matrix size
($m=25$e$+6$ and $n=50$).   The {\tt MvTransMv} operation requires a
global reduce, while the DGEMM operations do not require communication.  DCGS2
always outperforms CGS2 in these runs.  {\tt MvTimesMatAddMv} perform similarly
for both schemes on 96 and 192 nodes.  CGS2 exhibits an increase in speed that
matches DCGS2.  For a large number of rows or columns, DCGS2 obtains about $66$
GigaFlops/sec per node at 192 nodes or about 6\% of the single GPU sustained
execution rate of 200 GigaFlops/sec.

\ignore{
\begin{figure*}[htb]
        \centering
        \begin{subfigure}[b]{0.45\textwidth}
        \includegraphics[width=0.9\textwidth]{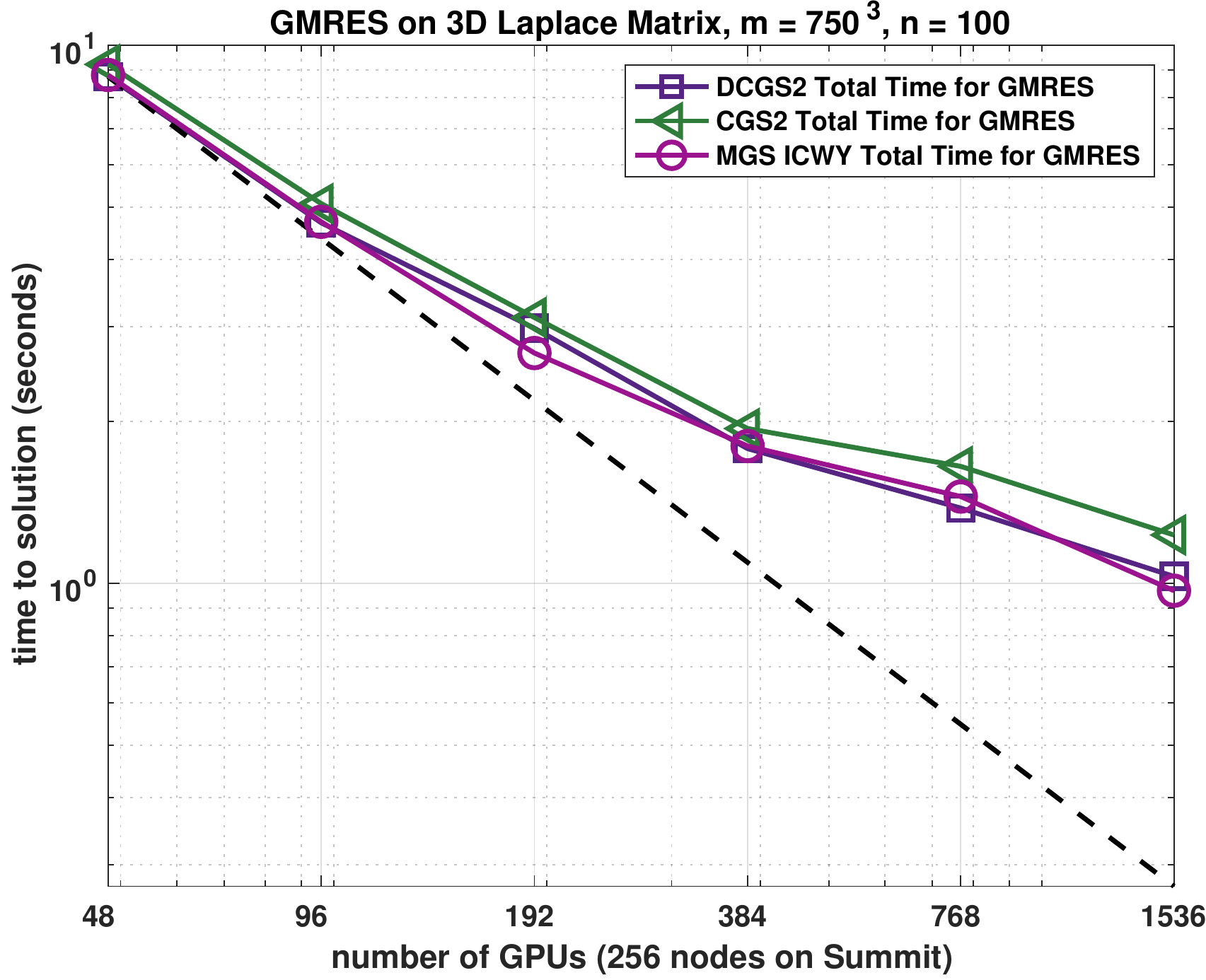}
        \caption{\label{fig:gmres_solve_time_750}GMRES up to 256 nodes (6 GPUs per node).}
        \end{subfigure}
        \quad\quad
        \begin{subfigure}[b]{0.45\textwidth}
        \includegraphics[width=0.9\textwidth]{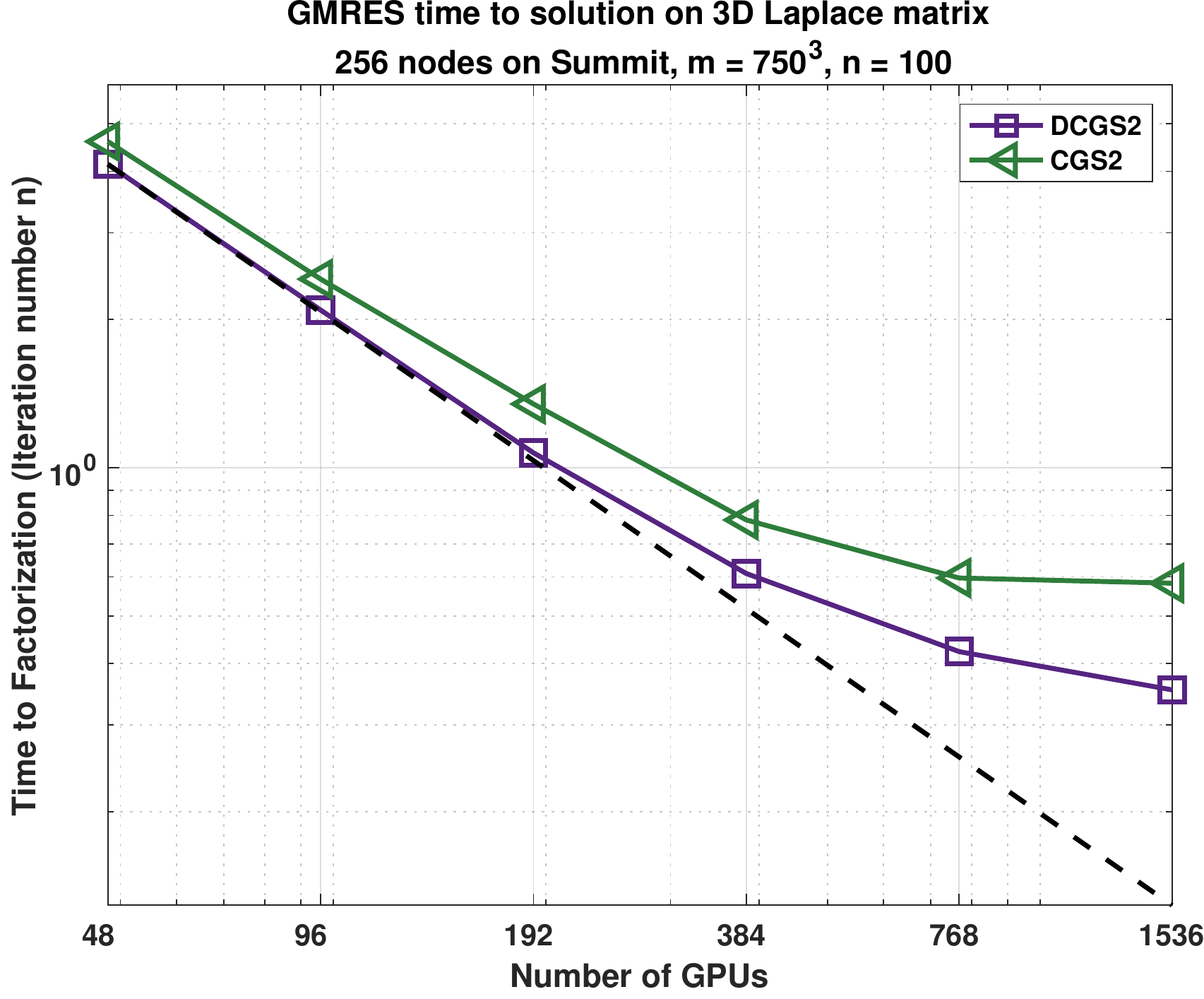}
        \caption{\label{fig:belos_solve_time_750}Belos up to 256 nodes (6 GPUs per node).}
        \end{subfigure}
        \caption{\label{fig:perf-time-belos-750}
         Time to solution of GMRES and Belos solves. $n=750^3$, $m=100$  (256 nodes).}
\end{figure*}

\begin{figure*}[htb]
        \centering
        \begin{subfigure}[b]{0.45\textwidth}
        \includegraphics[width=0.9\textwidth]{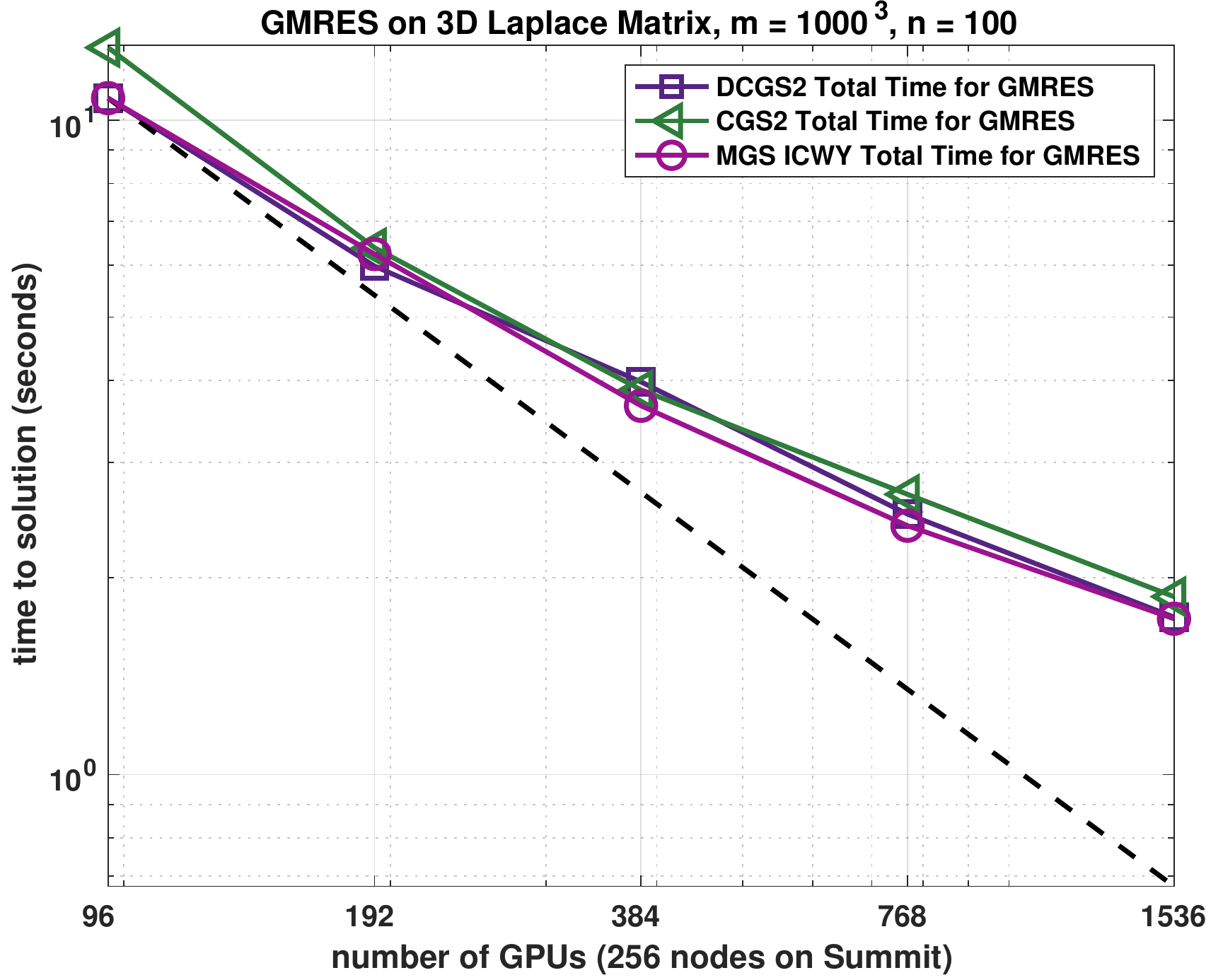}
        \caption{\label{fig:gmres_solve_time_1000}GMRES up to 256 nodes (6 GPUs per node).}
        \end{subfigure}
        \quad\quad
        \begin{subfigure}[b]{0.45\textwidth}
        \includegraphics[width=0.9\textwidth]{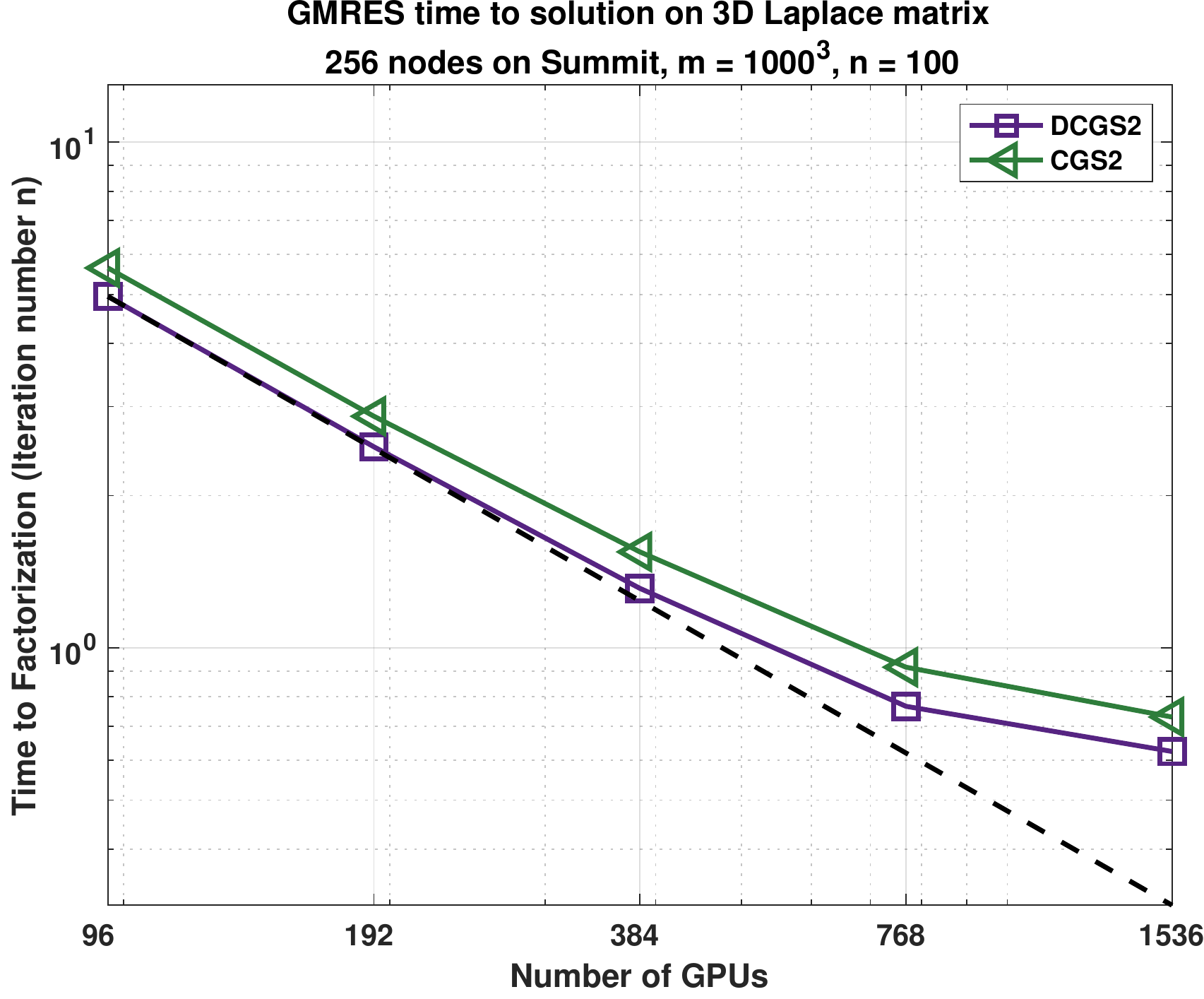}
        \caption{\label{fig:belos_solve_time_1000}Belos solve up to 256 nodes (6 GPUs per node).}
        \end{subfigure}
        \caption{\label{fig:perf-time-belos-1000} 
         Time to solution of GMRES and Belos solves. $n=1000^3$, $m=100$  (256 nodes).}
\end{figure*}
}
\begin{figure}[htb]
        \centering
        \includegraphics[width=0.45\textwidth]{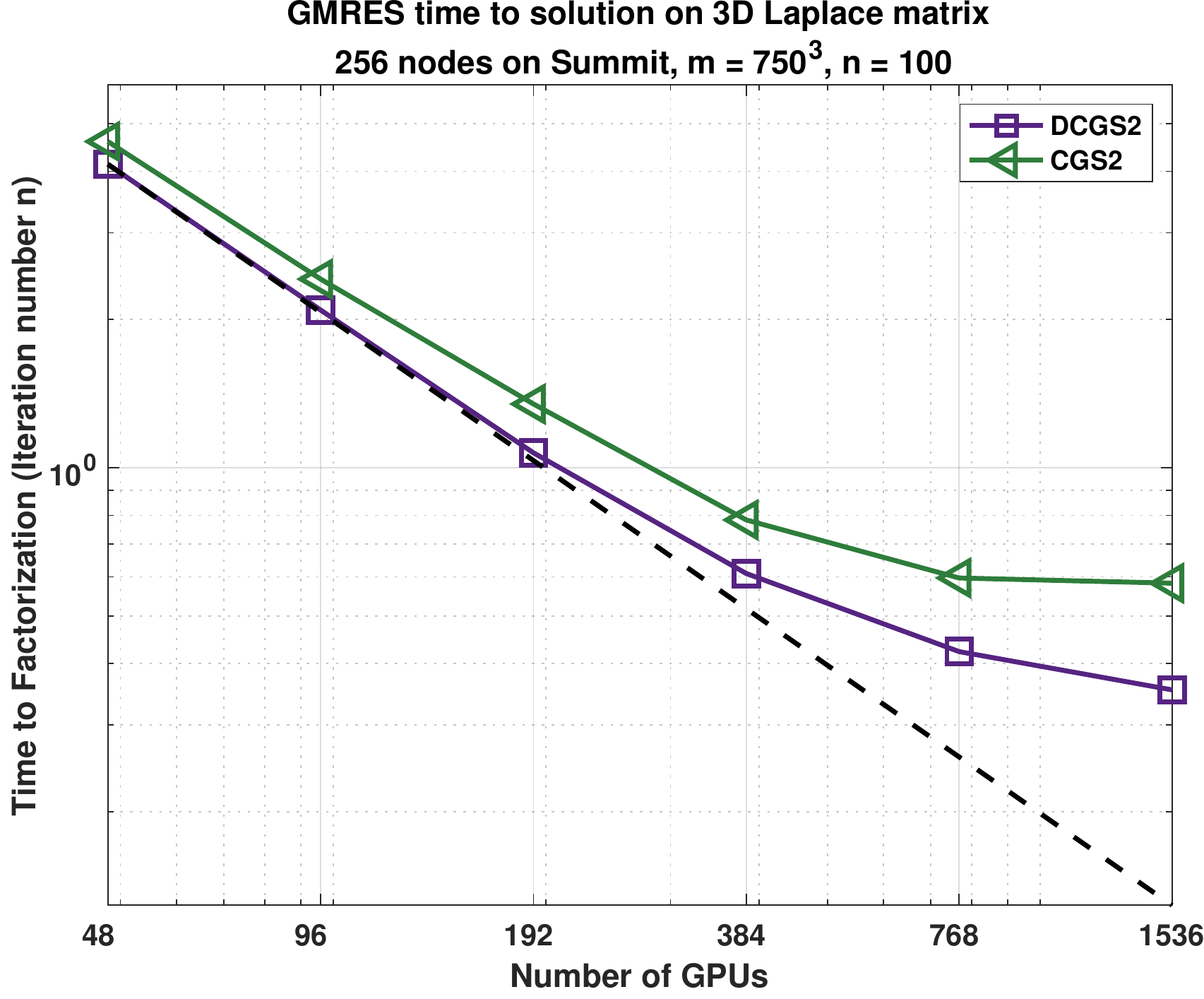}
        \caption{\label{fig:perf-time-belos-750}
         Time to solution of GMRES. $m=750^3$, $n=100$  (256 nodes, 6 GPUs per node).}
\end{figure}

\begin{figure}[htb]
        \centering
        \includegraphics[width=0.45\textwidth]{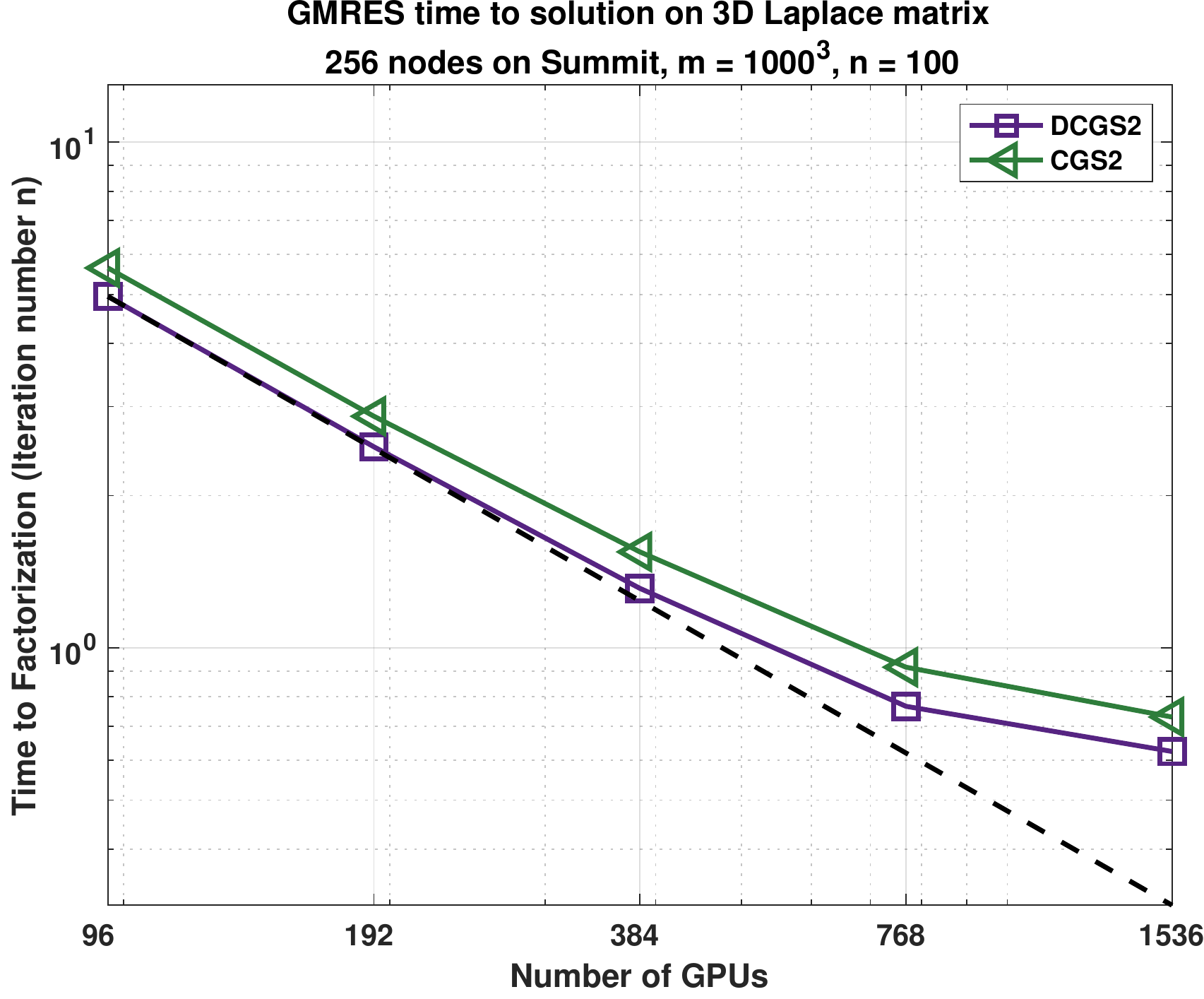}
        \caption{\label{fig:perf-time-belos-1000} 
         Time to solution of GMRES. $m=1000^3$, $n=100$  (256 nodes, 6 GPUs per node).}
\end{figure}

Figures~\ref{fig:perf-time-belos-750} and \ref{fig:perf-time-belos-1000} are strong scaling experiments for the 3D Laplace equation with dimension $m=750^3$ and $m=1000^3$. The simulations employ from 8 to 256 Summit compute nodes.  The GMRES solver is run for one trial for each of the node counts using 6 GPUs per node, in non-dedicated runs on Summit.
A fixed number of $n=100$ iterations are performed, without a restart,
and a preconditioner is not applied. Clearly, the DCGS2-GMRES yields
lower run times and exhibits better strong-scaling characteristics.

\section{Conclusion}
For distributed-memory computation, two passes of classical Gram-Schmidt (CGS2)
was the method of choice for Krylov solvers requiring machine precision level
representation errors and loss of orthogonality. However, the algorithm
requires three global reductions for each column of the $QR$ decomposition
computed and thus the strong-scaling behavior can deviate substantially from
linear as the number of MPI ranks increases on Exascale class supercomputers
such as the ORNL Summit.  In this paper, a new variant of CGS2 that requires
only one global reduction per iteration was applied to the Arnoldi-$QR$
algorithm.  Our numerical results have demonstrated that DCGS2 obtains the same
loss of orthogonality and representation error as CGS2, while our
strong-scaling studies on the Summit supercomputer demonstrate that DCGS2 obtains
a speedup of $2\times$ faster compute time on a single GPU, and an even larger
speedup on an increasing number of GPUs, reaching $2.2\times$ lower execution
times on 192 GPUs.  The impact of DCGS2 on the strong scaling of Krylov linear
system solvers is currently being explored, and a block variant is also being
implemented following the review article of Carson et
al.~\cite{2020_arXive_Carson}.
The software employed for this paper is available on GitHub.

\section*{Acknowledgement}
This research was supported by the Exascale Computing Project (17-SC-20-SC), 
a collaborative effort of the U.S. Department of Energy Office of Science 
and the National Nuclear Security Administration.  
The National Renewable Energy Laboratory is operated by Alliance for Sustainable
Energy, LLC, for the U.S. Department of Energy (DOE) under Contract No.
DE-AC36-08GO28308.  Sandia National Laboratories is a multimission laboratory
managed and operated by National Technology \& Engineering Solutions of Sandia,
LLC, a wholly owned subsidiary of Honeywell International Inc., for the U.S.
Department of Energy National Nuclear Security Administration under contract
DE-NA0003525.

A portion of this research used resources of the Oak Ridge Leadership Computing
Facility, that is a DOE Office of Science User Facility supported under
Contract DE-AC05-00OR22725 and using computational resources
sponsored by the Department of Energy's Office of Energy Efficiency and
Renewable Energy and located at the National Renewable Energy Laboratory.

Julien Langou was supported by NSF award \#1645514.

\bibliography{biblio___orthogonalization}

\end{document}